# Dimension and basis construction for $C^2$-smooth isogeometric spline spaces over bilinear-like $G^2$ two-patch parameterizations

Mario Kapl[a,b,*], Vito Vitrih[c,d]

[a] *Johann Radon Institute for Computational and Applied Mathematics, Austrian Academy of Sciences, Linz, Austria*
[b] *Dipartimento di Matematica "F.Casorati", Università degli Studi di Pavia, Italy*
[c] *IAM and FAMNIT, University of Primorska, Koper, Slovenia*
[d] *Institute of Mathematics, Physics and Mechanics, Ljubljana, Slovenia*

**Abstract**

A particular class of planar two-patch geometries, called bilinear-like $G^2$ two-patch geometries, is introduced. This class includes the subclass of all bilinear two-patch parameterizations and possesses similar connectivity functions along the patch interface. It is demonstrated that the class of bilinear-like $G^2$ two-patch parameterizations is much wider than the class of bilinear parameterizations and can approximate with good quality given generic two-patch parameterizations.

We investigate the space of $C^2$-smooth isogeometric functions over this specific class of two-patch geometries. The study is based on the equivalence of the $C^2$-smoothness of an isogeometric function and the $G^2$-smoothness of its graph surface (cf. [12, 20]). The dimension of the space is computed and an explicit basis construction is presented. The resulting basis functions possess simple closed form representations, have small local supports, and are well-conditioned. In addition, we introduce a subspace whose basis functions can be generated uniformly for all possible configurations of bilinear-like $G^2$ two-patch parameterizations. Numerical results obtained by performing $L^2$-approximation indicate that already the subspace possesses optimal approximation properties.

*Keywords:* isogeometric analysis, geometric continuity, second order continuity, two-patch domain, bilinear-like
*2000 MSC:* 65D17, 65N30, 68U07

## 1. Introduction

In Isogeometric Analysis (IgA), which was introduced in [14], partial differential equations (PDEs) are solved by using functions obtained from the parameterization (usually by polynomial or rational splines) of the computational domain [3, 9, 14]. This approach enables a direct link of the CAD representation of the computational domain with the numerical simulation of PDEs. High order PDEs can be solved by using their weak forms obtained by standard Galerkin discretization, see e.g. [2, 34]. This requires isogeometric spline spaces of $C^s$-smoothness, $s \geq 1$. While $C^1$-smooth isogeometric functions are needed to solve 4-th order PDEs (cf. [1, 4, 15, 20, 21, 22]), $C^2$-smooth functions are required for solving 6-th order PDEs by means of IgA (cf. [2, 11, 18, 19, 29, 34, 36]).

*Corresponding author
*Email addresses:* mario.kapl@ricam.oeaw.ac.at (Mario Kapl), vito.vitrih@upr.si (Vito Vitrih)

The study of such $C^s$-smooth isogeometric spline spaces on multi-patch domains is directly linked with the concept of geometric continuity (cf. [30, 31]). An isogeometric function is $C^s$-smooth on the given multi-patch domain if and only if its graph surface over the multi-patch domain is $G^s$-smooth (i.e., geometrically continuous of order $s$), see [12, 20]. Therefore, $C^s$-smooth isogeometric functions have been also called $C^s$-smooth geometrically continuous isogeometric functions [15, 18, 19, 20].

A widely used approach to construct $C^s$-smooth isogeometric functions on a given planar multi-patch geometry is to study the entire space of $C^s$-smooth isogeometric functions and to generate a basis of this space. The resulting basis is then used to describe the geometry of the given computational domain and to perform isogeometric simulation on this domain. In case of $s = 1$, it was shown in [8], that the class of so-called analysis-suitable $G^1$ multi-patch parameterizations is exactly the one which allows the construction of $C^1$-smooth isogeometric spaces with optimal approximation properties. Thereby, this class of parameterizations consists of particular $C^0$ multi-patch parameterizations and include amongst others the class of bilinear multi-patch geometries. Examples of constructions of $C^1$-smooth isogeometric spline spaces over specific analysis-suitable $G^1$ multi-patch parameterizations are [5, 8, 15, 17, 20, 25]. In case of $s = 2$, so far mainly bilinear multi-patch domains have been considered to generate $C^2$-smooth isogeometric spline spaces [18, 19]. Thereby, numerical experiments indicate that the constructed spaces possess again optimal approximation power.

Another approach to construct $C^s$-smooth isogeometric functions on planar multi-patch geometries, but not followed in this work, is to use domain parameterizations which are $C^1$-smooth along the patch interfaces except in the neighborhood of extraordinary vertices, see e.g. [26, 27, 28] for $s = 1$ and [35] for $s = 2$. In these examples the basis functions are derived from $G^s$-surface constructions which are able to handle extraordinary vertices.

In the present paper we focus on $s = 2$ and restrict ourselves to two-patch configurations. Already the two-patch case can handle some specific multi patch domains, see e.g. [19], and is an important pre-step to deal with multi-patch domains with extraordinary vertices, see e.g. [18]. More precisely, we consider a particular class of two-patch configurations, which are called bilinear-like $G^2$ two-patch parameterizations, and includes amongst others the bilinear two-patch parameterizations. The name is derived from the fact, that a bilinear-like $G^2$ parameterization possesses similar connectivity functions as a bilinear geometry along the patch interface.

The aim of this paper is to investigate the dimension of the space of $C^2$-smooth isogeometric functions over bilinear-like $G^2$ two-patch parameterizations and to provide an explicit and simple basis construction. The use of bilinear-like $G^2$ two-patch geometries instead of just bilinear two-patch parameterizations is advantageous, since it allows to model more general two-patch domains such as geometries with a curved interface and a curved boundary and still leads to $C^2$-smooth isogeometric spaces with optimal approximation properties. To demonstrate the potential and flexibility of the class of bilinear-like $G^2$ two-patch parameterizations to represent more general geometries, we describe a method for the construction of bilinear-like $G^2$ two-patch parameterizations approximating given generic two-patch geometries.

The present paper extends, generalizes and simplifies the work in [19], where the space of $C^2$-smooth biquintic and bisixtic isogeometric functions on bilinearly parameterized two-patch domains was considered, in several directions. First, as already mentioned above, our novel approach works for the much wider class of bilinear-like $G^2$ two-patch parameterizatons, which allows to deal with e.g. two-patch geometries with a curved interface and a curved boundary. Second, while the method here works for non-uniform splines of any bidegree $(p, p)$, $p \geq 5$, and of any regularity $r$, $2 \leq r < p-2$, within the single patches, the method in [19] is restricted so far to the case of uniform splines of bidegree $(p, p) = (5, 5)$ or $(p, p) = (6, 6)$ and regularity $r = 2$.

Third, in [19] a basis for the $C^2$-smooth isogeometric space is generated implicitly by using the concept of minimal determining sets (cf. [5, 23]) for the involved spline coefficients. Disadvantages of the resulting basis are the implicit representation, the required nullspace computation and a large support (usually over the whole interface) of the basis functions defined across the interface, which are not really desired properties for solving high order PDEs in IgA. In contrast, we present here an explict basis construction with closed form representations for the single basis functions. These functions have small, local supports, are well-conditioned and their spline coefficients can be simply computed by means of interpolation. This is achieved by using and developing similar tools as presented in [17] for the case of $C^1$-smooth isogeometric functions on analysis-suitable $G^1$ two-patch parameterizations. We use the constructed $C^2$-smooth isogeometric basis functions to perform $L^2$-approximation on different bilinear-like $G^2$ two-patch parameterizations. Analogous to the bilinear case [19], the numerical results indicate optimal approximation properties of the considered $C^2$-smooth isogeometric spline spaces. Moreover, we are able to introduce a particular subspace of the $C^2$-smooth space, whose basis functions have even a simpler and more uniform representation. We numerically show that already this subspace possesses optimal approximation order

Fourth, in [19] the study of the dimension of the $C^2$-smooth space is based on the categorization of four different possible configurations of bilinear two-patch geometries, and requires a lot of symbolic computation. Although we consider in the present paper the more general class of bilinear-like $G^2$ two-patch parameterizations, the study of the dimension is shortened and simplified, and works now uniformly for all possible configurations.

The remainder of the paper is organized as follows. Section 2 recalls the concept of $C^2$-smooth isogeometric functions over two-patch domains. In Section 3 and Section 4, we introduce the class of bilinear-like $G^2$ two-patch parameterizations, which is a generalization of the class of bilinear two-patch geometries studied in [19], and consider the associated $C^2$-smooth isogeometric spline spaces. A method to generate bilinear-like $G^2$ two-patch parameterizations from given generic two-patch parameterizations is described in Appendix A. While the dimension of the considered $C^2$-smooth isogeometric spline spaces is investigated in Section 5, an explicit basis construction is presented in Section 6. Furthermore, Section 7 outlines the approximation power of the constructed spaces. Finally, we conclude the paper in Section 8.

## 2. $C^2$ isogeometric spaces over two-patch parameterizations

*Spline spaces.* Throughout the paper, we will use a similar notation as introduced in [17]. Let $p \geq 5$, $2 \leq r \leq p-3$ and $k \in \mathbb{N}_0$. We define the open knot vector $\mathcal{T}_k^{p,r} = (t_0^{p,r}, t_1^{p,r}, \ldots, t_{2p+1+k(p-r)}^{p,r})$ as

$$\mathcal{T}_k^{p,r} = (\underbrace{0, 0, \ldots, 0}_{(p+1)-\text{times}}, \underbrace{\tau_1, \tau_1, \ldots, \tau_1}_{(p-r)-\text{times}}, \underbrace{\tau_2, \tau_2, \ldots, \tau_2}_{(p-r)-\text{times}}, \ldots, \underbrace{\tau_k, \tau_k, \ldots, \tau_k}_{(p-r)-\text{times}}, \underbrace{1, 1, \ldots, 1}_{(p+1)-\text{times}}),$$

where $0 < \tau_i < \tau_{i+1} < 1$ for all $1 \leq i \leq k-1$ and $k$ is the number of different inner knots. Let $\mathcal{S}(\mathcal{T}_k^{p,r}, [0,1])$ be the univariate spline space of degree $p$ defined on the interval $[0,1]$ with respect to the open knot vector $\mathcal{T}_k^{p,r}$, and let $N_i^{p,r}$, $i = 0, 1, \ldots, p + k(p-r)$, be the associated B-splines. Note that functions in the spline space $\mathcal{S}(\mathcal{T}_k^{p,r}, [0,1])$ are $C^r$-smooth at all inner knots.

In addition, we consider knot vectors, which are obtained from the knot vector $\mathcal{T}_k^{p,r}$ by increasing the multiplicity of single inner knots. For $\ell, \ell' \in \{1, 2, \ldots, k\}$ and $\ell \neq \ell'$, we obtain the knot vectors $\mathcal{T}_{k,\ell}^{p,r}$, $\mathcal{T}_{k,\ell,\ell}^{p,r}$ and $\mathcal{T}_{k,\ell,\ell'}^{p,r}$ by inserting the knot $\tau_\ell$ one more time, the knot $\tau_\ell$ two more times and each of the two knots $\tau_\ell$, $\tau_{\ell'}$ one more time, respectively. The resulting spline spaces $\mathcal{S}(\mathcal{T}_{k,\ell}^{p,r}, [0,1])$, $\mathcal{S}(\mathcal{T}_{k,\ell,\ell}^{p,r}, [0,1])$ and $\mathcal{S}(\mathcal{T}_{k,\ell,\ell'}^{p,r}, [0,1])$ are again $C^r$-smooth at all inner knots except at those where the multiplicity has been increased. At these knots the spline space is only $C^{r-1}$-smooth or $C^{r-2}$-smooth, depending if the corresponding knot has been inserted once or twice.

We denote by $\mathcal{S}(\mathcal{T}_k^{p,r}, [0,1]^2)$ the bivariate tensor-product spline space of bidegree $(p,p)$ on the unit square $[0,1]^2$ which is composed of the univariate spline space $\mathcal{S}(\mathcal{T}_k^{p,r}, [0,1])$, i.e.,

$$\mathcal{S}(\mathcal{T}_k^{p,r}, [0,1]^2) = \mathcal{S}(\mathcal{T}_k^{p,r}, [0,1]) \bigotimes \mathcal{S}(\mathcal{T}_k^{p,r}, [0,1]),$$

and possesses the tensor-product B-splines $N_{i,j}^{p,r} = N_i^{p,r} N_j^{p,r}$, $i, j = 0, 1, \ldots, p + k(p-r)$. Moreover, we denote by $\mathcal{P}^p([0,1])$ the univariate polynomial space of degree $p$ on the interval $[0,1]$ and by $\mathcal{P}^p([0,1]^2)$ the corresponding tensor-product polynomial space of bidegree $(p,p)$ on the unit square $[0,1]^2$.

*Two-patch parameterizations.* Let $\boldsymbol{F}^{(L)}$ and $\boldsymbol{F}^{(R)}$ be two regular and bijective geometry mappings

$$\boldsymbol{F}^{(S)} : [0,1]^2 \to \Omega^{(S)}, \quad \boldsymbol{F}^{(S)} \in \mathcal{S}(\mathcal{T}_k^{p,r}, [0,1]^2) \times \mathcal{S}(\mathcal{T}_k^{p,r}, [0,1]^2), \quad S \in \{L, R\},$$

with spline representations

$$\boldsymbol{F}^{(S)}(u,v) = \sum_{i=0}^{p+k(p-r)} \sum_{j=0}^{p+k(p-r)} \boldsymbol{c}_{i,j}^{(S)} N_{i,j}^{p,r}(u,v), \quad \boldsymbol{c}_{i,j}^{(S)} \in \mathbb{R}^2.$$

The two images $\boldsymbol{F}^{(S)}([0,1]^2)$, $S \in \{L, R\}$, describe quadrilateral patches $\Omega^{(S)} \subset \mathbb{R}^2$, i.e.,

$$\Omega^{(S)} = \boldsymbol{F}^{(S)}([0,1]^2),$$

and the union of these two patches defines a two-patch domain $\Omega \subset \mathbb{R}^2$, i.e.,

$$\Omega = \Omega^{(L)} \cup \Omega^{(R)}.$$

For the sake of simplicity, we assume that the two-patches share exactly one whole edge as the common interface $\Gamma = \Omega^{(L)} \cap \Omega^{(R)}$ (see Fig. 1) which is parameterized by

$$\boldsymbol{F}_0(v) = \boldsymbol{F}^{(L)}(0, v) = \boldsymbol{F}^{(R)}(0, v), \quad v \in [0,1].$$

In addition, we denote by $\boldsymbol{F}$ the two-patch parameterization consisting of the two geometry mappings $\boldsymbol{F}^{(L)}$ and $\boldsymbol{F}^{(R)}$.

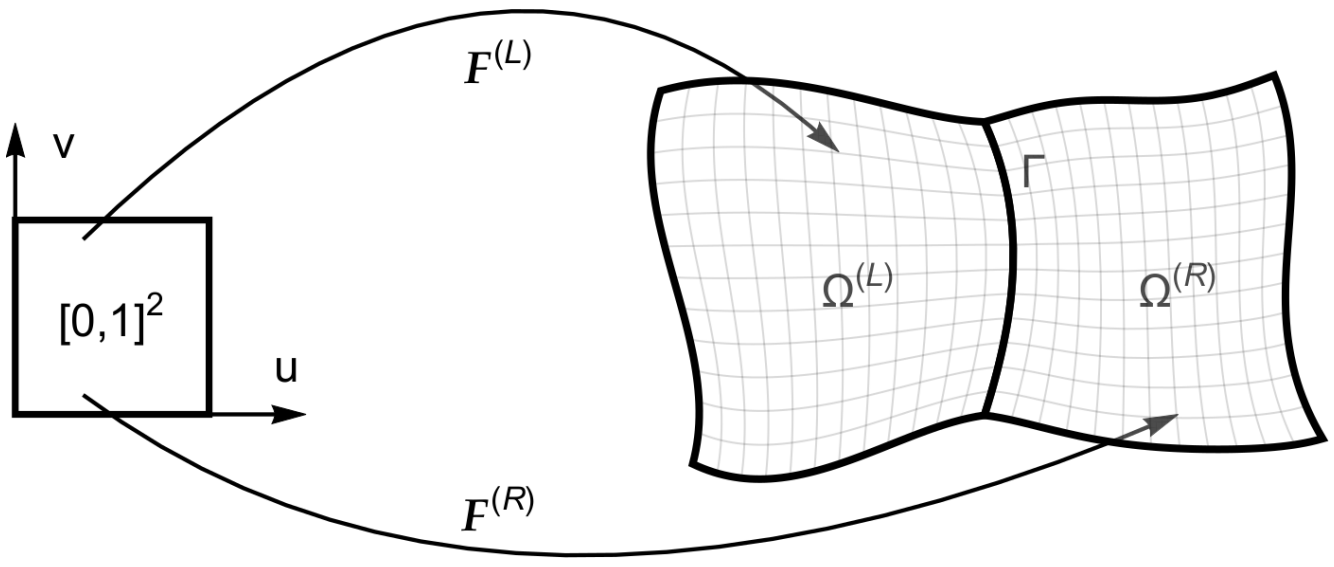

Figure 1: The two-patch domain $\Omega = \Omega^{(L)} \cup \Omega^{(R)}$ is described by geometry mappings $\boldsymbol{F}^{(L)}$ and $\boldsymbol{F}^{(R)}$.

*Isogeometric spaces.* Let $\mathcal{V}$ be the space of isogeometric functions on $\Omega$ with respect to the two-patch parameterization $\boldsymbol{F}$, i.e.,

$$\mathcal{V} = \{\phi : \Omega \to \mathbb{R} \mid \phi \circ \boldsymbol{F}^{(S)} \in \mathcal{S}(\mathcal{T}_k^{p,r}, [0,1]^2),\ S \in \{L, R\}\}.$$

Consider an isogeometric function $\phi \in \mathcal{V}$. For each patch $\Omega^{(S)}$, $S \in \{L, R\}$, the graph surface patch $\boldsymbol{\Sigma}^{(S)} : [0,1]^2 \to \Omega^{(S)} \times \mathbb{R}$ is given by

$$\boldsymbol{\Sigma}^{(S)}(u, v) = (\boldsymbol{F}^{(S)}(u, v), g^{(S)}(u, v))^T,$$

where $g^{(S)} = \phi \circ \boldsymbol{F}^{(S)} \in \mathcal{S}(\mathcal{T}_k^{p,r}, [0,1]^2)$, possessing the spline representation

$$g^{(S)}(u, v) = \sum_{i=0}^{p+k(p-r)} \sum_{j=0}^{p+k(p-r)} d_{i,j}^{(S)} N_{i,j}^{p,r}(u, v), \quad d_{i,j}^{(S)} \in \mathbb{R}.$$

Then, the graph surface $\boldsymbol{\Sigma} \subset \Omega \times \mathbb{R}$ of the isogemetric function $\phi$ is the surface, which is composed of the two graph surface patches $\boldsymbol{\Sigma}^{(L)}$ and $\boldsymbol{\Sigma}^{(R)}$.

*Geometric continuity conditions.* Let $\mathcal{V}^s$, $s \in \mathbb{N}_0$, be the subspace of $\mathcal{V}$ containing the globally $C^s$-smooth isogeometric functions $\phi \in \mathcal{V}$, i.e.,

$$\mathcal{V}^s = \mathcal{V} \cap C^s(\Omega).$$

The concept of geometric continuity (cf. [13, 31]) can be used to fully describe this space: *An isogeometric function $\phi \in \mathcal{V}$ belongs to the space $\mathcal{V}^s$ if and only if its graph surface $\mathbf{\Sigma}$ is $G^s$-smooth*, see e.g. [12, 20].

Below, we will focus on $s = 2$. To ensure that the graph surface $\mathbf{\Sigma}$ is $G^2$-smooth, the two graph surface patches $\mathbf{\Sigma}^{(L)}$ and $\mathbf{\Sigma}^{(R)}$ have to be joint along their common interface with $G^2$-continuity. This can be expressed by the following conditions (cf. [18, 19]):

$$\mathbf{\Sigma}^{(L)}(0,v) = \mathbf{\Sigma}^{(R)}(0,v), \quad v \in [0,1], \tag{1}$$

$$\det[D_u\mathbf{\Sigma}^{(L)}(0,v), D_u\mathbf{\Sigma}^{(R)}(0,v), D_v\mathbf{\Sigma}^{(L)}(0,v)] = 0, \quad v \in [0,1], \tag{2}$$

and

$$\det[\bar{\mathbf{\Xi}}(v), D_u\mathbf{\Sigma}^{(L)}(0,v), D_v\mathbf{\Sigma}^{(L)}(0,v)] = 0, \quad v \in [0,1], \tag{3}$$

where

$$\begin{aligned}\bar{\mathbf{\Xi}}(v) \;=\; & (\bar{\alpha}^{(L)}(v))^2 D_{uu}\mathbf{\Sigma}^{(R)}(0,v) - ((\bar{\alpha}^{(R)}(v))^2 D_{uu}\mathbf{\Sigma}^{(L)}(0,v) \\ & +2\bar{\alpha}^{(R)}(v)\bar{\beta}(v)D_{uv}\mathbf{\Sigma}^{(L)}(0,v) + \bar{\beta}^2(v)D_{vv}\mathbf{\Sigma}^{(L)}(0,v)),\end{aligned}$$

with

$$\bar{\alpha}^{(S)}(v) = \det[D_u\boldsymbol{F}^{(S)}(0,v), \boldsymbol{F}_0'(v)], \quad S \in \{L, R\},$$

and

$$\bar{\beta}(v) = \det[D_u\boldsymbol{F}^{(L)}(0,v), D_u\boldsymbol{F}^{(R)}(0,v)].$$

The conditions (1)-(3) are equivalent to the existence of functions $\alpha^{(L)}$, $\alpha^{(R)}$, $\beta$, $\widehat{\alpha}^{(L)}$, $\eta$, $\theta : [0,1] \to \mathbb{R}$ such that for all $v \in [0,1]$ the following conditions hold:

$$\alpha^{(L)}(v)\alpha^{(R)}(v) < 0, \tag{4}$$

$$\mathbf{\Sigma}^{(L)}(0,v) = \mathbf{\Sigma}^{(R)}(0,v),$$

$$\alpha^{(R)}(v)D_u\mathbf{\Sigma}^{(L)}(0,v) - \alpha^{(L)}(v)D_u\mathbf{\Sigma}^{(R)}(0,v) + \beta(v)D_v\mathbf{\Sigma}^{(L)}(0,v) = \mathbf{0},$$

and

$$\widehat{\alpha}^{(L)}(v)\mathbf{\Xi}(v) + \eta(v)D_u\mathbf{\Sigma}^{(L)}(0,v) + \theta(v)D_v\mathbf{\Sigma}^{(L)}(0,v) = \mathbf{0},$$

where

$$\begin{aligned}\mathbf{\Xi}(v) \;=\; & (\alpha^{(L)}(v))^2 D_{uu}\mathbf{\Sigma}^{(R)}(0,v) - ((\alpha^{(R)}(v))^2 D_{uu}\mathbf{\Sigma}^{(L)}(0,v) \\ & +2\alpha^{(R)}(v)\beta(v)D_{uv}\mathbf{\Sigma}^{(L)}(0,v) + \beta^2(v)D_{vv}\mathbf{\Sigma}^{(L)}(0,v)).\end{aligned}$$

Since the geometry mappings $\boldsymbol{F}^{(L)}$ and $\boldsymbol{F}^{(R)}$ are given, the functions $\alpha^{(L)}$, $\alpha^{(R)}$, $\beta$, $\widehat{\alpha}^{(L)}$, $\eta$ and $\theta$ possess the form

$$\alpha^{(L)}(v) = \gamma_1(v)\bar{\alpha}^{(L)}(v), \quad \alpha^{(R)}(v) = \gamma_1(v)\bar{\alpha}^{(R)}(v), \quad \beta(v) = \gamma_1(v)\bar{\beta}(v),$$

$$\widehat{\alpha}^{(L)}(v) = \gamma_2(v)\bar{\alpha}^{(L)}(v), \quad \eta(v) = \gamma_2(v)\det[\boldsymbol{F}_0'(v), \boldsymbol{Z}(v)]$$

and

$$\theta(v) = \gamma_2(v)\det[\boldsymbol{Z}(v), D_u\boldsymbol{F}^{(L)}(0,v)],$$

respectively, where $\gamma_i : [0,1] \to \mathbb{R}$, $i = 1,2$, is any function with $\gamma_i(v) \neq 0$, and

$$\begin{aligned} \boldsymbol{Z}(v) &= (\alpha^{(L)}(v))^2 D_{uu}\boldsymbol{F}^{(R)}(0,v) - ((\alpha^{(R)}(v))^2 D_{uu}\boldsymbol{F}^{(L)}(0,v) \\ &\quad +2\alpha^{(R)}(v)\beta(v)D_{uv}\boldsymbol{F}^{(L)}(0,v) + \beta^2(v)\boldsymbol{F}_0''(v)). \end{aligned} \tag{5}$$

Moreover, the function $\beta$ can be written as

$$\beta(v) = \alpha^{(L)}(v)\beta^{(R)}(v) - \alpha^{(R)}(v)\beta^{(L)}(v), \tag{6}$$

where $\beta^{(L)}$, $\beta^{(R)} : [0,1] \to \mathbb{R}$ are non-unique functions, see e.g. [8, 31]. One possible choice of the functions $\beta^{(S)}$, $S \in \{L, R\}$, is given by

$$\beta^{(S)}(v) = \frac{D_u\boldsymbol{F}^{(S)}(0,v) \cdot \boldsymbol{F}_0'(v)}{||\boldsymbol{F}_0'(v)||^2}. \tag{7}$$

Furthermore, the selected geometry mappings $\boldsymbol{F}^{(L)}$ and $\boldsymbol{F}^{(R)}$ ensure that condition (4) is fulfilled, which also implies that for $v \in [0,1]$ $\alpha^{(S)}(v) \neq 0$, $S \in \{L, R\}$. Finally, we obtain that an isogeometric function $\phi \in \mathcal{V}$ is $C^2$-smooth on $\Omega$ (i.e. $\phi \in \mathcal{V}^2$) if and only if

$$g^{(L)}(0,v) = g^{(R)}(0,v), \tag{8}$$

$$\alpha^{(R)}(v)D_u g^{(L)}(0,v) - \alpha^{(L)}(v)D_u g^{(R)}(0,v) + \beta(v)D_v g^{(L)}(0,v) = 0, \tag{9}$$

and

$$\widehat{\alpha}^{(L)}(v)w(v) + \eta(v)D_u g^{(L)}(0,v) + \theta(v)D_v g^{(L)}(0,v) = 0, \tag{10}$$

where

$$\begin{aligned} w(v) &= (\alpha^{(L)}(v))^2 D_{uu}g^{(R)}(0,v) - ((\alpha^{(R)}(v))^2 D_{uu}g^{(L)}(0,v) \\ &\quad + 2\alpha^{(R)}(v)\beta(v)D_{uv}g^{(L)}(0,v) + \beta^2(v)D_{vv}g^{(L)}(0,v)). \end{aligned} \tag{11}$$

Condition (8) guarantees that $\phi$ is $C^0$-smooth, condition (9) additionally ensures that $\phi$ is $C^1$-smooth, and condition (10) finally implies that $\phi$ is $C^2$-smooth.

In case of bilinear geometry mappings $\boldsymbol{F}^{(L)}$ and $\boldsymbol{F}^{(R)}$, i.e. $\boldsymbol{F}^{(L)}, \boldsymbol{F}^{(R)} \in \mathcal{P}^1([0,1]^2)$), the functions $\widehat{\alpha}^{(L)}$, $\eta$ and $\theta$ can be expressed by means of linear functions $\alpha^{(L)}$, $\alpha^{(R)}$, $\beta^{(L)}$ and $\beta^{(R)}$.

**Proposition 1.** *Let $\boldsymbol{F}^{(L)}, \boldsymbol{F}^{(R)} \in \mathcal{P}^1([0,1]^2)$, $\gamma_1(v) = \gamma_2(v) = 1$ and $\beta^{(L)}$ and $\beta^{(R)}$ as given in* (7)*. Then, we obtain that*

$$\alpha^{(L)}, \alpha^{(R)}, \beta^{(L)}, \beta^{(R)} \in \mathcal{P}^1([0,1]),$$

*and*

$$\widehat{\alpha}^{(L)}(v) = \alpha^{(L)}(v), \tag{12}$$

$$\eta(v) = 2(\alpha^{(L)})'(v)\alpha^{(R)}(v)\beta(v), \tag{13}$$

$$\theta(v) = 2\left(\alpha^{(L)}(v)(\beta^{(L)})'(v) - (\alpha^{(L)})'(v)\beta^{(L)}(v)\right)\alpha^{(R)}(v)\beta(v). \tag{14}$$

*Proof.* This can be shown by a straightforward computation. □

## 3. Bilinear-like $G^2$ two-patch parameterizations

Below, we will consider a particular class of two-patch parameterizations, which possesses functions $\alpha^{(L)}$, $\alpha^{(R)}$, $\beta$, $\widehat{\alpha}^{(L)}$, $\eta$ and $\theta$ with a similar form as for bilinear two-patch parameterizations. Let us define this class of two-patch parameterizations.

**Definition 1** (Bilinear-like $G^2$ two-patch parameterization)**.** Two-patch parameterization $\boldsymbol{F}$ consisting of geometry mappings $\boldsymbol{F}^{(L)}$, $\boldsymbol{F}^{(R)} \in \mathcal{S}(\mathcal{T}_k^{p,r},[0,1]^2) \times \mathcal{S}(\mathcal{T}_k^{p,r},[0,1]^2)$ is called *bilinear-like* $G^2$ if there exist functions $\alpha^{(L)}, \alpha^{(R)}$, $\beta^{(L)}$, $\beta^{(R)} \in \mathcal{P}^1([0,1])$ satisfying equations

$$\boldsymbol{F}^{(L)}(0,v) = \boldsymbol{F}^{(R)}(0,v), \tag{15}$$

$$\alpha^{(R)}(v) D_u \boldsymbol{F}^{(L)}(0,v) - \alpha^{(L)}(v) D_u \boldsymbol{F}^{(R)}(0,v) + \beta(v) D_v \boldsymbol{F}^{(L)}(0,v) = \boldsymbol{0}, \tag{16}$$

$$\widehat{\alpha}^{(L)}(v)\boldsymbol{Z}(v) + \eta(v) D_u \boldsymbol{F}^{(L)}(0,v) + \theta(v) D_v \boldsymbol{F}^{(L)}(0,v) = \boldsymbol{0}, \tag{17}$$

for $v \in [0,1]$, where $\boldsymbol{Z}$, $\beta$, $\widehat{\alpha}^{(L)}$, $\eta$ and $\theta$ are given by (5), (6), (12), (13) and (14), respectively.

In order to obtain $C^2$-smooth isogeometric spaces with optimal approximation properties, we restrict the functions $\alpha^{(L)}, \alpha^{(R)}$, $\beta^{(L)}$, $\beta^{(R)}$ to be linear, which is motivated by the study of analysis suitable $G^1$ parameterizations for the case of $C^1$-smooth isogeometric spaces, cf. [8].

Two instances of bilinear and bilinear-like $G^2$ two-patch parameterizations with the same functions $\alpha^{(L)}$, $\alpha^{(R)}$, $\beta$, $\widehat{\alpha}^{(L)}$, $\eta$ and $\theta$ are given in Fig. A.4. An advantage of using bilinear-like $G^2$ two-patch parameterizations instead of bilinear two-patch parameterization is the possibility to represent two-patch domains with a curved interface and a curved boundary. A simple method for the construction of bilinear-like $G^2$ two-patch parameterizations is presented in Appendix A, and allows to generate bilinear-like $G^2$ two-patch geometries which approximate with good quality given generic two-patch parameterizations.

**Lemma 1.** *Let $\boldsymbol{F}$ be a bilinear-like $G^2$ two-patch parameterization consisting of geometry mappings $\boldsymbol{F}^{(L)}, \boldsymbol{F}^{(R)} \in \mathcal{S}(\mathcal{T}_k^{p,r},[0,1]^2)\times\mathcal{S}(\mathcal{T}_k^{p,r},[0,1]^2)$ with functions $\alpha^{(L)}, \alpha^{(R)}$, $\beta^{(L)}$, $\beta^{(R)} \in \mathcal{P}^1([0,1])$. An isogeometric function $\phi \in \mathcal{V}$ belongs to the space $\mathcal{V}^2$ if and only if*

$$g^{(L)}(0,v) = g^{(R)}(0,v), \tag{18}$$

$$\frac{D_u g^{(L)}(0,v) - \beta^{(L)}(v) g_0'(v)}{\alpha^{(L)}(v)} = \frac{D_u g^{(R)}(0,v) - \beta^{(R)}(v) g_0'(v)}{\alpha^{(R)}(v)}, \tag{19}$$

*and*

$$\frac{D_{uu} g^{(L)}(0,v) - (\beta^{(L)}(v))^2 g_0''(v) - 2\alpha^{(L)}(v)\beta^{(L)}(v) g_1'(v)}{(\alpha^{(L)}(v))^2} = \frac{D_{uu} g^{(R)}(0,v) - (\beta^{(R)}(v))^2 g_0''(v) - 2\alpha^{(R)}(v)\beta^{(R)}(v) g_1'(v)}{(\alpha^{(R)}(v))^2} \tag{20}$$

*for all* $v \in [0,1]$*, where* $g_0, g_1, g_2 : [0,1] \to \mathbb{R}$ *are (possibly rational) functions representing the two equally valued terms in* (18)*,* (19) *and* (20)*, i.e.,*

$$g_0(v) = g^{(S)}(0,v), \quad g_1(v) = \frac{D_u g^{(S)}(0,v) - \beta^{(S)}(v) g_0'(v)}{\alpha^{(S)}(v)}, \quad S \in \{L,R\},$$

$$g_2(v) = \frac{D_{uu} g^{(S)}(0,v) - (\beta^{(S)}(v))^2 g_0''(v) - 2\alpha^{(S)}(v)\beta^{(S)}(v) g_1'(v)}{(\alpha^{(S)}(v))^2}, \quad S \in \{L,R\}.$$

*Proof.* Since the equations (8) and (18) are equal and the equivalence of (9) and (19) has already been shown in [8], it remains to show that the equations (10) and (20) are equivalent. To keep the formulas short, we will omit the arguments of the involved functions.

Using the fact that the functions $w$, $\widehat{\alpha}^{(L)}$, $\eta$ and $\theta$ are given by (11), (12), (13) and (14), respectively, the equation (10) is equivalent to

$$\begin{aligned}&(\alpha^{(L)})^3 D_{uu} g^{(R)} - \alpha^{(L)}(\alpha^{(R)})^2 D_{uu} g^{(L)} - 2\alpha^{(L)}\alpha^{(R)}\beta D_{uv} g^{(L)} - \alpha^{(L)}\beta^2 D_{vv} g^{(L)} \\ &+2(\alpha^{(L)})'\alpha^{(R)}\beta D_u g^{(L)} + 2(\alpha^{(L)}(\beta^{(L)})' - (\alpha^{(L)})'\beta^{(L)})\alpha^{(R)}\beta D_v g^{(L)} = 0.\end{aligned}$$

Since we have

$$D_v g^{(L)} = g_0', \quad D_u g^{(L)} = \beta^{(L)} g_0' + \alpha^{(L)} g_1,$$

and therefore

$$D_{vv} g^{(L)} = g_0'', \quad D_{uv} g^{(L)} = (\beta^{(L)})' g_0' + \beta^{(L)} g_0'' + (\alpha^{(L)})' g_1 + \alpha^{(L)} g_1',$$

we obtain that

$$\begin{aligned}&(\alpha^{(L)})^3 D_{uu} g^{(R)} - \alpha^{(L)}(\alpha^{(R)})^2 D_{uu} g^{(L)} - 2\alpha^{(L)}\alpha^{(R)}\beta((\beta^{(L)})' g_0' + \beta^{(L)} g_0'' + (\alpha^{(L)})' g_1 + \alpha^{(L)} g_1') \\ &- \alpha^{(L)}\beta^2 g_0'' + 2(\alpha^{(L)})'\alpha^{(R)}\beta(\beta^{(L)} g_0' + \alpha^{(L)} g_1) + 2(\alpha^{(L)}(\beta^{(L)})' - (\alpha^{(L)})'\beta^{(L)})\alpha^{(R)}\beta g_0' = 0,\end{aligned}$$

which simplifies to

$$(\alpha^{(L)})^3 D_{uu} g^{(R)} - \alpha^{(L)}(\alpha^{(R)})^2 D_{uu} g^{(L)} - 2\alpha^{(L)}\alpha^{(R)}\beta^{(L)}\beta g_0'' - 2(\alpha^{(L)})^2\alpha^{(R)}\beta g_1' - \alpha^{(L)}\beta^2 g_0'' = 0. \tag{21}$$

By dividing equation (21) by $\alpha^{(L)}$ and reordering the expression, we obtain

$$(\alpha^{(L)})^2 D_{uu} g^{(R)} = (\alpha^{(R)})^2 D_{uu} g^{(L)} + 2\alpha^{(R)}\beta^{(L)}\beta g_0'' + 2\alpha^{(L)}\alpha^{(R)}\beta g_1' + \beta^2 g_0''. \tag{22}$$

By substituting the function $\beta$ by (6), the equation (22) is equivalent to

$$\begin{aligned}(\alpha^{(L)})^2 D_{uu} g^{(R)} &= (\alpha^{(R)})^2 D_{uu} g^{(L)} + 2\alpha^{(R)}\beta^{(L)}(\alpha^{(L)}\beta^{(R)} - \alpha^{(R)}\beta^{(L)}) g_0'' \\ &\quad +2\alpha^{(L)}\alpha^{(R)}(\alpha^{(L)}\beta^{(R)} - \alpha^{(R)}\beta^{(L)}) g_1' + (\alpha^{(L)}\beta^{(R)} - \alpha^{(R)}\beta^{(L)})^2 g_0'',\end{aligned}$$

which simplifies to

$$\begin{aligned}(\alpha^{(L)})^2 D_{uu} g^{(R)} - (\alpha^{(L)})^2(\beta^{(R)})^2 g_0'' - 2(\alpha^{(L)})^2\alpha^{(R)}\beta^{(R)} g_1' &= \\ (\alpha^{(R)})^2 D_{uu} g^{(L)} - (\alpha^{(R)})^2(\beta^{(L)})^2 g_0'' - 2(\alpha^{(R)})^2\alpha^{(L)}\beta^{(L)} g_1'&.\end{aligned} \tag{23}$$

By diving equation (23) by $(\alpha^{(L)}\alpha^{(R)})^2$, we obtain equation (20). □

**Remark 1.** Lemma 1 can be extended to all two-patch parameterizations $\boldsymbol{F}$ for which there exist any functions $\alpha^{(L)}$, $\alpha^{(R)}$, $\beta^{(L)}$ and $\beta^{(R)}$ such that (5), (6) and (12)-(17) hold.

In the following of the paper we will restrict ourselves to bilinear-like $G^2$ two-patch geometries $\boldsymbol{F}$ with geometry mappings $\boldsymbol{F}^{(L)}, \boldsymbol{F}^{(R)} \in \mathcal{S}(\mathcal{T}_k^{p,r}, [0,1]^2) \times \mathcal{S}(\mathcal{T}_k^{p,r}, [0,1]^2)$ and functions $\alpha^{(L)}, \alpha^{(R)}$, $\beta^{(L)}$, $\beta^{(R)} \in \mathcal{P}^1([0,1])$, and we will investigate the space $\mathcal{V}^2$ for these particular class of two-patch geometries.

## 4. $C^2$ isogeometric spaces over bilinear-like $G^2$ two-patch parameterizations

We will investigate the space $\mathcal{V}^2$ for bilinear-like $G^2$ two-patch parameterizations and will decompose $\mathcal{V}^2$ into the direct sum of two subspaces. Before, we introduce the two index sets

$$\mathcal{I} = \{(S,i,j) \mid S \in \{L,R\} \text{ and } i,j = 0,1,\ldots,p+k(p-r)\},$$

and

$$\mathcal{I}_\Gamma = \{(S,i,j) \mid S \in \{L,R\},\ i = 0,1,2 \text{ and } j = 0,1,\ldots,p+k(p-r)\}.$$

Thereby, the set $\mathcal{I}$ collects the indices of all spline control points $d_{i,j}^{(S)}$, $S \in \{L,R\}$, (of a function $\phi \in \mathcal{V}^2$) and the set $\mathcal{I}_\Gamma$ collects the indices of all spline control points $d_{i,j}^{(S)}$, $S \in \{L,R\}$, which are influenced by the interface $\Gamma$, i.e., appear in the equations (18)-(20). The space $\mathcal{V}^2$ can be decomposed into

$$\mathcal{V}^2 = \mathcal{V}_1^2 \oplus \mathcal{V}_2^2, \tag{24}$$

where the subspaces $\mathcal{V}_1^2$ and $\mathcal{V}_2^2$ are given by

$$\mathcal{V}_1^2 = \{\phi \in \mathcal{V}^2 \mid d_{i,j}^{(S)} = 0 \text{ for } (S,i,j) \in \mathcal{I}_\Gamma\},$$

and

$$\mathcal{V}_2^2 = \{\phi \in \mathcal{V}^2 \mid d_{i,j}^{(S)} = 0 \text{ for } (S,i,j) \in \mathcal{I} \setminus \mathcal{I}_\Gamma\},$$

respectively. The decomposition (24) implies that

$$\dim \mathcal{V}^2 = \dim \mathcal{V}_1^2 + \dim \mathcal{V}_2^2.$$

**Remark 2.** We also call the functions of $\mathcal{V}_1^2$ and $\mathcal{V}_2^2$ functions of the first and second kind, respectively (cf. [19]).

*The space* $\mathcal{V}_1^2$. Since the functions in $\mathcal{V}_1^2$ are not affected by the interface $\Gamma$, the space $\mathcal{V}_1^2$ consists of all functions $\phi \in \mathcal{V}^2$ with the property that the functions $g^{(S)} = \phi \circ \boldsymbol{F}^{(S)}$, $S \in \{L,R\}$, possess a spline representation

$$g^{(S)}(u,v) = \sum_{i=3}^{p+k(p-r)} \sum_{j=0}^{p+k(p-r)} d_{i,j}^{(S)} N_{i,j}^{p,r}(u,v), \quad d_{i,j}^{(S)} \in \mathbb{R}.$$

This implies that the dimension of $\mathcal{V}_1^2$ is equal to the cardinality of the set $\mathcal{I} \setminus \mathcal{I}_\Gamma$ , i.e., $\dim \mathcal{V}_1^2 = |\mathcal{I} \setminus \mathcal{I}_\Gamma|$, which gives us the following lemma.

**Lemma 2.** *The dimension of* $\mathcal{V}_1^2$ *is equal to*

$$\dim \mathcal{V}_1^2 = 2(p-2+k(p-r))(p+1+k(p-r)).$$

A possible basis of the space $\mathcal{V}_1^2$ is given by the collection of the isogeometric functions $\phi_{i,j}^{(L)}$ and $\phi_{i,j}^{(R)}$, $i=3,4,\ldots,p+k(p-r)$, $j=0,1,\ldots,p+k(p-r)$, which are determined by

$$\phi_{i,j}^{(L)}:\ g^{(L)}(u,v) = N_{i,j}^{p,r}(u,v),\ g^{(R)}(u,v)=0,$$

and

$$\phi_{i,j}^{(R)}:\ g^{(L)}(u,v) = 0,\ g^{(R)}(u,v)=N_{i,j}^{p,r}(u,v),$$

respectively.

*The space* $\mathcal{V}_2^2$. We will present an explicit representation of the function $\phi \in \mathcal{V}_2^2$. Let us first introduce some notation, which will be needed in the rest of the paper. Let $q$ be the polynomial $q = \gcd(\alpha^{(L)}, \alpha^{(R)})$ with the leading coefficient being equal to one, and let

$$\widetilde{\alpha}^{(S)} = \frac{\alpha^{(S)}}{q}, \quad S \in \{L,R\}. \tag{25}$$

Furthermore, let us define the polynomial $h$ as

$$h = \begin{cases} 1, & \text{if } \gcd(\beta^{(L)},\beta^{(R)}) = \gcd(\alpha^{(L)},\alpha^{(R)}) = q, \\ q, & \text{otherwise.} \end{cases} \tag{26}$$

Additionally, let $d_\alpha = \max\left(\deg(\alpha^{(L)}), \deg(\alpha^{(R)})\right)$, $d_{\widetilde{\alpha}} = \max\left(\deg(\widetilde{\alpha}^{(L)}), \deg(\widetilde{\alpha}^{(R)})\right)$ and $d_h = \deg h$. We will also need later the index sets

$$\mathcal{Z}_\beta = \{i \in \{1,2,\ldots,k\},\ \beta(\tau_i)=0\}, \quad \mathcal{Z}_\beta^C := \{1,2,\ldots,k\}\backslash \mathcal{Z}_\beta,$$

and assume that these sets are ordered. Let us denote the cardinality of the set $\mathcal{Z}_\beta$ by $z_\beta$, i.e., $z_\beta = |\mathcal{Z}_\beta|$. Since for bilinear-like $G^2$ two-patch parameterizations $\beta \in \mathcal{P}^2([0,1])$, $z_\beta \in \{0,1,2,k\}$. Depending on $\beta$, we define the knot vector $\widetilde{\mathcal{T}}_k^p$ as

$$\begin{cases} \mathcal{T}_k^{p,r} & \text{if } \beta = 0, \\ \mathcal{T}_k^{p,r+1} & \text{if } \beta \neq 0 \text{ and } z_\beta = 0, \\ \mathcal{T}_{k,\ell}^{p,r+1} & \text{if } \beta \neq 0 \text{ and } z_\beta = 1, \\ \mathcal{T}_{k,\ell,\ell'}^{p,r+1} & \text{if } \beta \neq 0 \text{ and } z_\beta = 2, \end{cases}$$

where $\ell, \ell' \in \{1,2,\ldots,k\}$ with $\ell \neq \ell'$ are the indices of possible roots $\tau_\ell$ and $\tau_{\ell'}$ of $\beta$.

Recall that Lemma 1 states the $C^2$ conditions for $\phi \in \mathcal{V}_2^2$ with respect to the two associated functions $g^{(S)} = \phi \circ \boldsymbol{F}^{(S)}$, $S \in \{L,R\}$. We use these conditions to define functions $\widetilde{g}_0$, $\widetilde{g}_1$, $\widetilde{g}_2 : [0,1] \to \mathbb{R}$. Let $g_0$ be defined as in Lemma 1 and let $\widetilde{g}_0 = g_0$, i.e.,

$$g^{(L)}(0,v) = g^{(R)}(0,v) = g_0(v) = \widetilde{g}_0(v). \tag{27}$$

Using notation (25), we can rewrite (19), and define $\widetilde{g}_1$ as

$$\frac{D_u g^{(L)}(0,v) - \beta^{(L)}(v)\widetilde{g}_0{}'(v)}{q(v)\,\widetilde{\alpha}^{(L)}(v)} = \frac{D_u g^{(R)}(0,v) - \beta^{(R)}(v)\widetilde{g}_0{}'(v)}{q(v)\,\widetilde{\alpha}^{(R)}(v)} = g_1(v) = \frac{\widetilde{g}_1(v)}{q(v)}. \tag{28}$$

Therefrom, we obtain that

$$D_u g^{(S)}(0,v) = \alpha^{(S)}(v) g_1(v) + \beta^{(S)}(v) g_0'(v) = \widetilde{\alpha}^{(S)}(v)\widetilde{g}_1(v) + \beta^{(S)}(v)\widetilde{g}_0{}'(v),\ S \in \{L,R\}. \tag{29}$$

Further, the function $\widetilde{g}_2$ is defined by reformulating equation (20) in the way that

$$\begin{aligned}
&\frac{D_{uu} g^{(L)}(0,v) - \left(\beta^{(L)}(v)\right)^2 \widetilde{g}_0{}''(v) - 2\widetilde{\alpha}^{(L)}(v)\beta^{(L)}(v)\left(\widetilde{g}_1{}'(v) - \widetilde{g}_1(v)\frac{q'(v)}{q(v)}\right)}{q^2(v)\left(\widetilde{\alpha}^{(L)}(v)\right)^2} \\
&= \frac{D_{uu} g^{(R)}(0,v) - \left(\beta^{(R)}(v)\right)^2 \widetilde{g}_0{}''(v) - 2\widetilde{\alpha}^{(R)}(v)\beta^{(R)}(v)\left(\widetilde{g}_1{}'(v) - \widetilde{g}_1(v)\frac{q'(v)}{q(v)}\right)}{q^2(v)\left(\widetilde{\alpha}^{(R)}(v)\right)^2} \qquad (30) \\
&= g_2(v) = \frac{\widetilde{g}_2(v)}{q^2(v)}. \qquad (31)
\end{aligned}$$

This implies for $S \in \{L,R\}$ that

$$D_{uu} g^{(S)}(0,v) = \left(\beta^{(S)}(v)\right)^2 \widetilde{g}_0{}''(v) + 2\widetilde{\alpha}^{(S)}(v)\beta^{(S)}(v)\left(\widetilde{g}_1{}'(v) - \widetilde{g}_1(v)\frac{q'(v)}{q(v)}\right) + \left(\widetilde{\alpha}^{(S)}(v)\right)^2 \widetilde{g}_2(v). \tag{32}$$

We can now specify an explicit representation of $\phi \in \mathcal{V}_2^2$.

**Lemma 3.** *For a function $\phi \in \mathcal{V}_2^2$, the associated functions $g^{(S)} = \phi \circ \boldsymbol{F}^{(S)}$, $S \in \{L,R\}$, can be represented as*

$$\begin{aligned}
g^{(S)}(u,v) = \widetilde{g}_0(v) M_0(u) + \left(\widetilde{\alpha}^{(S)}(v)\widetilde{g}_1(v) + \beta^{(S)}(v)\widetilde{g}_0{}'(v)\right) M_1(u) + \Big(\left(\beta^{(S)}(v)\right)^2 \widetilde{g}_0{}''(v) + \\
2\widetilde{\alpha}^{(S)}(v)\beta^{(S)}(v)\big(\widetilde{g}_1{}'(v) - \widetilde{g}_1(v)\frac{q'(v)}{q(v)}\big) + \left(\widetilde{\alpha}^{(S)}(v)\right)^2 \widetilde{g}_2(v)\Big) M_2(u),
\end{aligned} \tag{33}$$

*where $\widetilde{g}_0$, $\widetilde{g}_1$ and $\widetilde{g}_2$ are defined in* (27), (28) *and* (31)*, respectively, and*

$$M_0(u) = \sum_{i=0}^{2} N_i^{p,r}(u),\ M_1(u) = \frac{\tau_1}{p}\left(N_1^{p,r}(u) + 2N_2^{p,r}(u)\right),\ M_2(u) = \frac{\tau_1^2}{p(p-1)} N_2^{p,r}(u). \tag{34}$$

*Proof.* By means of Taylor expansion of $g^{(S)}(u,v)$ at $(u,v) = (0,v)$ and using equations

(29) and (32), we obtain that

$$
\begin{aligned}
g^{(S)}(u,v) &= g^{(S)}(0,v) + D_u g^{(S)}(0,v)u + D_{uu} g^{(S)}(0,v)\frac{u^2}{2} + \mathcal{O}(u^3)\\
&= \widetilde{g}_0(v) + \left(\widetilde{\alpha}^{(S)}(v)\widetilde{g}_1(v) + \beta^{(S)}(v)\widetilde{g}_0{}'(v)\right)u + \left(\left(\beta^{(S)}(v)\right)^2 \widetilde{g}_0{}''(v) + 2\widetilde{\alpha}^{(S)}(v)\beta^{(S)}(v)\cdot\right.\\
&\quad \left.\left(\widetilde{g}_1{}'(v) - \widetilde{g}_1(v)\frac{q'(v)}{q(v)}\right) + \left(\widetilde{\alpha}^{(S)}(v)\right)^2 \widetilde{g}_2(v)\right)\frac{u^2}{2} + \mathcal{O}(u^3)\\
&= \widetilde{g}_0(v)M_0(u) + \left(\widetilde{\alpha}^{(S)}(v)\widetilde{g}_1(v) + \beta^{(S)}(v)\widetilde{g}_0{}'(v)\right)M_1(u) + \left(\left(\beta^{(S)}(v)\right)^2 \widetilde{g}_0{}''(v) + \right.\\
&\quad \left. 2\widetilde{\alpha}^{(S)}(v)\beta^{(S)}(v)\left(\widetilde{g}_1{}'(v) - \widetilde{g}_1(v)\frac{q'(v)}{q(v)}\right) + \left(\widetilde{\alpha}^{(S)}(v)\right)^2 \widetilde{g}_2(v)\right)M_2(u),
\end{aligned}
$$

where

$$
M_j(u) = \sum_{i=0}^{2} \lambda_{j,i} N_i^{p,r}(u), \quad j = 0,1,2.
$$

We used the fact that $g^{(S)}(u,v)$, $S \in \{L,R\}$, are splines with only three columns of nonzero spline coefficients $d_{i,j}^{(S)}$. Unknown parameters $\lambda_{j,i}$ are then determined by the conditions

$$
M_j(0) = \delta_{j,0}, \quad M_j'(0) = \delta_{j,1}, \quad M_j''(0) = \delta_{j,2},
$$

where $\delta_{j,\ell}$ is the Kronecker delta, and the B-spline properties

$$
N_i^{p,r}(0) = \begin{cases} 1, & i = 0,\\ 0, & i > 0, \end{cases} \quad \left(N_i^{p,r}\right)'(0) = \begin{cases} \frac{-p}{\tau_1}, & i = 0,\\ \frac{p}{\tau_1}, & i = 1,\\ 0, & i > 1, \end{cases} \quad \left(N_i^{p,r}\right)''(0) = \begin{cases} \frac{p(p-1)}{\tau_1^2}, & i = 0,\\ \frac{-2p(p-1)}{\tau_1^2}, & i = 1,\\ \frac{p(p-1)}{\tau_1^2}, & i = 2,\\ 0, & i > 2. \end{cases}
$$

□

Applying relations

$$
g_i = \frac{\widetilde{g}_i}{q^i}, \quad i = 0,1,2,
$$

in (33) we get the following corollary.

**Corollary 1.** *For a function $\phi \in \mathcal{V}_2^2$, the associated functions $g^{(S)} = \phi \circ \boldsymbol{F}^{(S)}$, $S \in \{L,R\}$, can be represented as*

$$
\begin{aligned}
g^{(S)}(u,v) = g_0(v)M_0(u) + \left(\alpha^{(S)}(v)g_1(v) + \beta^{(S)}(v)g_0'(v)\right)M_1(u) + &\\
\left(\left(\beta^{(S)}(v)\right)^2 g_0''(v) + 2\alpha^{(S)}(v)\beta^{(S)}(v)g_1'(v) + \left(\alpha^{(S)}(v)\right)^2 g_2(v)\right)M_2(u),&
\end{aligned} \tag{35}
$$

*where $g_0$, $g_1$ and $g_2$ are defined as in Lemma 1.*

**Remark 3.** In contrast to $\widetilde{g}_1$ and $\widetilde{g}_2$, the functions $g_1$ and $g_2$ can be rational.

Note that vice versa, not any arbitrary triplet of functions $(\widetilde{g}_0, \widetilde{g}_1, \widetilde{g}_2)$ in (33) would provide a function $g^{(S)} = \phi \circ \boldsymbol{F}^{(S)} \in \mathcal{S}(\mathcal{T}_k^{p,r}, [0,1]^2)$, and hence an isogeometric function $\phi \in \mathcal{V}_2^2$. Let $\widehat{\mathcal{V}_2^2}$ be the space of all possible triplets $(\widetilde{g}_0, \widetilde{g}_1, \widetilde{g}_2)$ such that $g^{(S)} \in \mathcal{S}(\mathcal{T}_k^{p,r}, [0,1]^2)$, $S \in \{L, R\}$, i.e.,

$$\widehat{\mathcal{V}_2^2} = \{(\widetilde{g}_0, \widetilde{g}_1, \widetilde{g}_2) \mid \text{both representations (33) of } g^{(S)},\ S \in \{L, R\}, \text{ belong to } \mathcal{S}(\mathcal{T}_k^{p,r}, [0,1]^2)\}.$$

We obtain the following theorem and corollary.

**Theorem 1.** *There is a one-to-one correspondence between the spaces $\mathcal{V}_2^2$ and $\widehat{\mathcal{V}_2^2}$.*

*Proof.* For any $\phi \in \mathcal{V}_2^2$ the functions $\widetilde{g}_0$, $\widetilde{g}_1$ and $\widetilde{g}_2$ are uniquely specified by (27), (28) and (31), respectively. Vice versa, for any $(\widetilde{g}_0, \widetilde{g}_1, \widetilde{g}_2) \in \widehat{\mathcal{V}_2^2}$ there exists a unique $\phi \in \mathcal{V}_2^2$, since representation (33) is uniquely determined by any possible choice of functions $\widetilde{g}_0$, $\widetilde{g}_1$ and $\widetilde{g}_2$. □

**Corollary 2.** *The dimension of $\mathcal{V}_2^2$ is equal to the dimension of $\widehat{\mathcal{V}_2^2}$.*

## 5. Dimension of the space $\mathcal{V}^2$ for bilinear-like $G^2$ two-patch geometries

The dimension of $\mathcal{V}_1^2$ has been already specified in Lemma 2. The dimension of $\mathcal{V}_2^2$ will be obtained by means of Corolloray 2. To compute the dimension of $\widehat{\mathcal{V}_2^2}$ we will consider the following three spaces

$$\begin{aligned}
\Gamma_0 &= \{\widetilde{g}_0 \mid (\widetilde{g}_0, \widetilde{g}_1, \widetilde{g}_2) \in \widehat{\mathcal{V}_2^2} \text{ for some } \widetilde{g}_1 \text{ and } \widetilde{g}_2\},\\
\Gamma_1 &= \{\widetilde{g}_1 \mid (0, \widetilde{g}_1, \widetilde{g}_2) \in \widehat{\mathcal{V}_2^2} \text{ for some } \widetilde{g}_2\},\\
\Gamma_2 &= \{\widetilde{g}_2 \mid (0, 0, \widetilde{g}_2) \in \widehat{\mathcal{V}_2^2}\},
\end{aligned}$$

and finally use the fact that

$$\dim \mathcal{V}_2^2 = \dim \widehat{\mathcal{V}_2^2} = \dim \Gamma_0 + \dim \Gamma_1 + \dim \Gamma_2. \tag{36}$$

In the next three lemmas we prove dimensions of spaces $\Gamma_0$, $\Gamma_1$ and $\Gamma_2$.

**Lemma 4.** *The space $\Gamma_2$ and its dimension are equal to*

$$\Gamma_2 = \mathcal{S}(\mathcal{T}_k^{p-2d_{\widetilde{\alpha}},r}, [0,1]) \quad \text{and} \quad \dim \Gamma_2 = k(p - 2d_{\widetilde{\alpha}} - r) + p + 1 - 2d_{\widetilde{\alpha}}.$$

*Proof.* By selecting $\widetilde{g}_0 = \widetilde{g}_1 = 0$ in (33) it follows that $\widetilde{g}_2 \in \mathcal{S}(\mathcal{T}_k^{p-2d_{\widetilde{\alpha}},r}, [0,1])$. The dimension formula follows directly from the definition of the spline space. □

**Lemma 5.** *It holds that*

$$\Gamma_1 = \{h\bar{g}_1 \mid \bar{g}_1 \in \mathcal{S}(\widetilde{\mathcal{T}_k}^{p-d_{\widetilde{\alpha}}-d_h}, [0,1])\},$$

*where $h$ is defined in* (26)*, and*

$$\dim \Gamma_1 = k(p - d_{\widetilde{\alpha}} - d_h - r - 1) + p - d_{\widetilde{\alpha}} - d_h + z_\beta + 1.$$

*Proof.* By selecting $\widetilde{g}_0 = 0$ and since $D_u g^{(S)}(0,v) \in \mathcal{S}\left(\mathcal{T}_k^{p,r}, [0,1]\right)$, (28) implies that (see [17, Lemma 4])

$$\widetilde{g}_1(v) = \frac{D_u g^{(S)}(0,v)}{\widetilde{\alpha}^{(S)}(v)} \in \mathcal{S}(\mathcal{T}_k^{p-d_{\widetilde{\alpha}},r}, [0,1]). \tag{37}$$

Moreover, representation (33) simplifies to

$$g^{(S)}(u,v) = \widetilde{\alpha}^{(S)}(v)\widetilde{g}_1(v)M_1(u) + \left(2\widetilde{\alpha}^{(S)}\beta^{(S)}\big(\widetilde{g}_1{}' - \widetilde{g}_1\frac{q'}{q}\big) + \big(\widetilde{\alpha}^{(S)}\big)^2\,\widetilde{g}_2\right)(v)\,M_2(u).$$

By (30) we obtain

$$(\widetilde{\alpha}^{(R)})^2(v)D_{uu}g^{(L)}(0,v) - (\widetilde{\alpha}^{(L)})^2(v)D_{uu}g^{(R)}(0,v) + 2\widetilde{\alpha}^{(R)}(v)\widetilde{\alpha}^{(L)}(v)\left(\widetilde{\alpha}^{(L)}(v)\beta^{(R)}(v)\right.$$
$$\left. - \widetilde{\alpha}^{(R)}(v)\beta^{(L)}(v)\right)\left(\widetilde{g}_1{}'(v) - \frac{\widetilde{g}_1(v)q'(v)}{q(v)}\right) = 0.$$

Therefrom, clearly

$$\widetilde{g}_1 = h\,\bar{g}_1, \quad \bar{g}_1 : [0,1] \to \mathbb{R}.$$

Since $D_{uu}g^{(S)}(0,v) \in \mathcal{S}\left(\mathcal{T}_k^{p,r}, [0,1]\right)$ and by using (37) it is straightforward to see, by considering all cases regarding $\beta$, that $\bar{g}_1 \in \mathcal{S}(\widetilde{\mathcal{T}_k}^{p-d_{\widetilde{\alpha}}-d_h}, [0,1])$ and therefore $\Gamma_1 \subseteq \{h\bar{g}_1 \mid \bar{g}_1 \in \mathcal{S}(\widetilde{\mathcal{T}_k}^{p-d_{\widetilde{\alpha}}-d_h}, [0,1])\}$. It remains to prove that $\Gamma_1 \supseteq \{h\bar{g}_1 \mid \bar{g}_1 \in \mathcal{S}(\widetilde{\mathcal{T}_k}^{p-d_{\widetilde{\alpha}}-d_h}, [0,1])\}$. More precisely, for every $\bar{g}_1 \in \mathcal{S}(\widetilde{\mathcal{T}_k}^{p-d_{\widetilde{\alpha}}-d_h}, [0,1])$ we have to find some $\widetilde{g}_2$, such that $(0, \widetilde{g}_1, \widetilde{g}_2) \in \widehat{\mathcal{V}_2^2}$. Let first $z_\beta = 0$. Then we can simply choose $\widetilde{g}_2 = 0$. Let now $z_\beta \neq 0$, which implies $|\mathcal{Z}_\beta| \neq 0$. Then we can represent $\bar{g}_1$ as

$$\bar{g}_1 = \sum_{i \in \mathcal{Z}_\beta} \bar{g}_{1,i}, \quad \bar{g}_{1,i} \in \mathcal{S}(\mathcal{T}_{k,i}^{p-d_{\widetilde{\alpha}}-d_h, r+1}, [0,1]),$$

and the appropriate choice for $\widetilde{g}_2$ is

$$\widetilde{g}_2(v) = -2h(v)\sum_{i \in \mathcal{Z}_\beta} \frac{\beta^{(L)}(\tau_i)}{\widetilde{\alpha}^{(L)}(\tau_i)} \bar{g}'_{1,i}(v).$$

The dimension formula for $\Gamma_1$ follows directly from the definition of the spline space. □

**Lemma 6.** *The space $\Gamma_0$ and its dimension are equal to*

$$\Gamma_0 = \mathcal{S}(\widetilde{\mathcal{T}}_k^p, [0,1]) \quad \textit{and} \quad \dim \Gamma_0 = k(p - r - 1) + p + z_\beta + 1.$$

*Proof.* By (33) clearly $\widetilde{g_0} \in \mathcal{S}(\mathcal{T}_k^{p,r},[0,1])$. Since $D_u g^{(S)}(0,v) \in \mathcal{S}(\mathcal{T}_k^{p,r},[0,1])$, relation (9) additionally implies $\widetilde{g_0} \in \mathcal{S}(\widetilde{\mathcal{T}}_k^{p},[0,1])$ and therefore $\Gamma_0 \subseteq \mathcal{S}(\widetilde{\mathcal{T}}_k^{p},[0,1])$. To prove $\Gamma_0 \supseteq \mathcal{S}(\widetilde{\mathcal{T}}_k^{p},[0,1])$, we have to find for every $\widetilde{g_0} \in \mathcal{S}(\widetilde{\mathcal{T}}_k^{p},[0,1])$ some $\widetilde{g_1}$ and $\widetilde{g_2}$, such that $(\widetilde{g_0},\widetilde{g_1},\widetilde{g_2}) \in \widehat{\mathcal{V}_2^2}$. We can represent $\widetilde{g_0}$ as

$$\widetilde{g_0} = \sum_{i \in \mathcal{Z}_\beta^C} g_{0,i} + \sum_{i \in \mathcal{Z}_\beta} \hat{g}_{0,i},$$

with $g_{0,i} \in \mathcal{S}(\mathcal{T}_{k,i}^{p,r+2},[0,1])$ and $\hat{g}_{0,i} \in \mathcal{S}(\mathcal{T}_{k,i,i}^{p,r+2},[0,1])$. The appropriate choices for $\widetilde{g_1}$ and $\widetilde{g_2}$ are then

$$\widetilde{g_1}(v) = -q(v)\left(\sum_{i \in \mathcal{Z}_\beta^C} \frac{\widetilde{\alpha}^{(R)}(\tau_i)\beta^{(L)}(\tau_i) + \widetilde{\alpha}^{(L)}(\tau_i)\beta^{(R)}(\tau_i)}{2\widetilde{\alpha}^{(R)}(\tau_i)\widetilde{\alpha}^{(L)}(\tau_i)q(\tau_i)} g'_{0,i}(v) + \sum_{i \in \mathcal{Z}_\beta} \frac{\beta^{(L)}(\tau_i)}{q(\tau_i)\widetilde{\alpha}^{(L)}(\tau_i)} \hat{g}'_{0,i}(v)\right),$$

and

$$\widetilde{g_2}(v) = \sum_{i \in \mathcal{Z}_\beta^C} \frac{\beta^{(L)}(\tau_i)\beta^{(R)}(\tau_i)}{\widetilde{\alpha}^{(L)}(\tau_i)\widetilde{\alpha}^{(R)}(\tau_i)} g''_{0,i}(v) + \sum_{i \in \mathcal{Z}_\beta} \frac{(\beta^{(L)}(\tau_i))^2}{(\widetilde{\alpha}^{(L)}(\tau_i))^2} \hat{g}''_{0,i}(v).$$

The dimension formula follows directly from the definition of the spline space. □

**Corollary 3.** *The dimension of $\mathcal{V}_2^2$ is equal to*

$$\dim \mathcal{V}_2^2 = (k+1)\,(3(p+1) - 3d_{\widetilde{\alpha}} - d_h) - (3r+5)k + 2z_\beta. \tag{38}$$

*Proof.* It follows directly by (36) and Lemma (4)-(6). □

Finally, we obtain the dimension of $\mathcal{V}^2$.

**Theorem 2.** *It holds that*

$$\dim \mathcal{V}^2 = \underbrace{2(p-2+k(p-r))(p+1+k(p-r))}_{\dim \mathcal{V}_1^2} + \underbrace{(k+1)\,(3(p+1) - 3d_{\widetilde{\alpha}} - d_h) - (3r+5)k + 2z_\beta}_{\dim \mathcal{V}_2^2}.$$

**Remark 4.** In [19] the dimension problem has been already considered for $r = 2$ and $p \in \{5, 6\}$ in case of bilinear two-patch geometries. Thereby, four different types of configurations (denoted as A-D) of the bilinear geometry mappings were used to compute the dimension. Table 1 shows the evaluation of our dimension formula of $\mathcal{V}_2^2$ for the four configurations and confirms that our results are in agreement with the ones in [19].

We present the dimension of $\mathcal{V}_2^2$ for two concrete instances of bilinear-like $G^2$ two-patch parameterizations.

**Example 1.** We consider the two bilinear-like $G^2$ two-patch parameterizations $\boldsymbol{F}$ visualized in Fig. A.4 (third row). Both geometries are generated from initial two-patch parameterizations $\widetilde{\boldsymbol{F}}$, see Fig. A.4 (second row), by using the least-squares approach. The detailed construction of both bilinear-like $G^2$ two-patch parameterizations as well as the exact representations of the involved geometries are given in Appendix A.

| Type of configuration | $d_{\widetilde{\alpha}}$ | $d_h$ | $\dim \mathcal{V}_2^2$ |
|---|---|---|---|
| Configuration A | 1 | 0 | $(k+1)3p - 11k + 2z_\beta$ |
| Configuration B | 0 | 0 | $(k+1)(3p+3) - 11k + 2z_\beta$ |
| Configuration C | 0 | 0 | $(k+1)(3p+3) - 11k + 2z_\beta$ |
| Configuration D | 0 | 1 | $(k+1)(3p+2) - 11k + 2z_\beta$ |

Table 1: Evaluation of formula (38) for $r = 2$ and the four different types of configurations of bilinear two-patch domains considered in [19].

The bilinear-like $G^2$ two-patch geometry (a) possesses linear functions $\alpha^{(L)}$, $\alpha^{(R)}$, $\beta^{(L)}$ and $\beta^{(R)}$, which are given by

$$\alpha^{(L)}(v) = -9-v, \;\; \alpha^{(R)}(v) = -\frac{3}{2}(-7+v), \;\; \beta^{(L)}(v) = \frac{1}{18}(-3+5v), \;\; \beta^{(R)}(v) = \frac{1}{12}(-1+3v).$$

This results in $q = h = 1$ and $\beta = \frac{1}{6}(15 - 32v + v^2)$, which implies that $d_{\widetilde{\alpha}} = 1$, $d_h = 0$ and $z_\beta = 0$. Therefore, the dimension of $\mathcal{V}_2^2$ is given by

$$\dim \mathcal{V}_2^2 = (k+1)(3(p+1)-3) - (3r+5)k = 3p + k(3p - 3r - 5).$$

In case of the bilinear-like $G^2$ two-patch parameterization (b), the linear functions $\alpha^{(L)}$, $\alpha^{(R)}$, $\beta^{(L)}$ and $\beta^{(R)}$ are given by

$$\alpha^{(L)}(v) = 9(-2+v), \quad \alpha^{(R)}(v) = -9(-2+v), \quad \beta^{(L)}(v) = \beta^{(R)}(v) = 1 - \frac{v}{2}.$$

It follows that $q = -2 + v$, $h = 1$ and $\beta = -9(-2+v)^2$ and hence $d_{\widetilde{\alpha}} = 0$, $d_h = 0$ and $z_\beta = 0$. This leads to

$$\dim \mathcal{V}_2^2 = 3(k+1)(p+1) - (3r+5)k = 3p + 3 + k(3p - 3r - 2).$$

## 6. Basis of the space $\mathcal{V}^2$ for bilinear-like $G^2$ two-patch geometries

A basis of the space $\mathcal{V}_1^2$ has been already generated in Section 4. In this section, we will present an explicit basis representation for the space $\mathcal{V}_2^2$.

*Explicit representation.* For a given bilinear-like $G^2$ two-patch geometry $\boldsymbol{F}$ with the geometry mappings $\boldsymbol{F}^{(L)}$ and $\boldsymbol{F}^{(R)}$ and functions $\alpha^{(L)}$, $\alpha^{(R)}$, $\beta^{(L)}, \beta^{(R)} \in \mathcal{P}^1([0,1])$, we will construct basis functions of $\mathcal{V}_2^2$ and provide explicit representations for these functions. For this purpose, particular sets of B-splines will be needed. For an ordered set $\mathcal{J} = \{j_0, j_1, \ldots\} \subseteq \{1, 2, \ldots, k\}$, we denote by $\mathcal{N}(\mathcal{T}_k^{p,r}, \mathcal{J})$ the set of $|\mathcal{J}|$ B-splines $\{N_0, N_1, \ldots, N_{|\mathcal{J}|-1}\}$, where $N_i$ is a selected B-spline $N_i \in \mathcal{S}(\mathcal{T}_{k,j_i}^{p,r}, [0,1])$ with the property $N_i(\tau_{j_i}) \neq 0$. Similarly, we define for $\mathcal{J}$ the set $\widetilde{\mathcal{N}}(\mathcal{T}_k^{p,r}, \mathcal{J})$ of $|\mathcal{J}|$ B-splines $\{\widetilde{N}_0, \widetilde{N}_1, \ldots, \widetilde{N}_{|\mathcal{J}|-1}\}$, where $\widetilde{N}_i$ is a selected B-spline $\widetilde{N}_i \in \mathcal{S}(\mathcal{T}_{k,j_i,j_i}^{p,r}, [0,1])$ with the property $\widetilde{N}_i(\tau_{j_i}) \neq 0$.

A basis of $\mathcal{V}_2^2$ is constructed as follows. We choose suitable function triplets $(\widetilde{g}_0, \widetilde{g}_1, \widetilde{g}_2)$, which determine via relation (33) basis functions of $\mathcal{V}_2^2$. Thereby, these function triplets are selected via basis for the spaces $\Gamma_i$, $i = 0,1,2$. For each $i = 0,1,2$, we denote by $\phi_{i,j}$, $j = 0,1,\ldots,\dim\Gamma_i - 1$, the corresponding basis functions of $\mathcal{V}_2^2$ and by $g_{i,j}^{(S)}$ the resulting spline functions $\phi_{i,j} \circ \boldsymbol{F}^{(S)}$, $S \in \{L,R\}$.

- $\Gamma_0$:

$$\{\phi_{0,j}\}_{j=0}^{n_0-1}: \ (\widetilde{g}_0, \widetilde{g}_1, \widetilde{g}_2) = \left(N_j^{p,r+2}, 0, 0\right), \quad n_0 = \dim(\mathcal{S}(\mathcal{T}_k^{p,r+2}, [0,1])),$$

$$\{\phi_{0,n_0+j}\}_{j=0}^{k-1}: \ (\widetilde{g}_0, \ \widetilde{g}_1, \widetilde{g}_2) = \left(N_j, -\frac{\widetilde{\alpha}^{(R)}(\tau_{j+1})\beta^{(L)}(\tau_{j+1}) + \widetilde{\alpha}^{(L)}(\tau_{j+1})\beta^{(R)}(\tau_{j+1})}{2\widetilde{\alpha}^{(R)}(\tau_{j+1})\widetilde{\alpha}^{(L)}(\tau_{j+1})q(\tau_{j+1})} \, q\, N_j', \right.$$
$$\left. \frac{\beta^{(L)}(\tau_{j+1})\beta^{(R)}(\tau_{j+1})}{\widetilde{\alpha}^{(L)}(\tau_{j+1})\widetilde{\alpha}^{(R)}(\tau_{j+1})} N_j'' \right), \quad N_j \in \mathcal{N}(\mathcal{T}_k^{p,r+2}, \{1,2,\ldots,k\}),$$

$$\{\phi_{0,n_0+k+j}\}_{j=0}^{z_\beta-1}: \ (\widetilde{g}_0, \widetilde{g}_1, \widetilde{g}_2) = \left(\widetilde{N}_j, -\frac{\beta^{(L)}(\tau_{i_j})}{q(\tau_{i_j})\widetilde{\alpha}^{(L)}(\tau_{i_j})} \, q\, \widetilde{N}_j', \ \frac{(\beta^{(L)}(\tau_{i_j}))^2}{(\widetilde{\alpha}^{(L)}(\tau_{i_j}))^2} \widetilde{N}_j'' \right),$$
$$\widetilde{N}_j \in \widetilde{\mathcal{N}}(\mathcal{T}_k^{p,r+2}, \mathcal{Z}_\beta), \quad \mathcal{Z}_\beta = \{i_0, i_1, \ldots, i_{z_\beta-1}\}.$$

- $\Gamma_1$:

$$\{\phi_{1,j}\}_{j=0}^{n_1-1}: \ (0, \widetilde{g}_1, \widetilde{g}_2) = \left(0, h\, N_j^{p-d_{\widetilde{\alpha}}-d_h, r+1}, 0\right), \quad n_1 = \dim(\mathcal{S}(\mathcal{T}_k^{p-d_{\widetilde{\alpha}}-d_h, r+1}, [0,1])),$$

$$\{\phi_{1,n_1+j}\}_{j=0}^{z_\beta-1}: (0, \widetilde{g}_1, \widetilde{g}_2) = \left(0, h\, N_j, -\frac{2\,\beta^{(L)}(\tau_{i_j})}{\widetilde{\alpha}^{(L)}(\tau_{i_j})} \, h\, N_j'\right), \ N_j \in \mathcal{N}(\mathcal{T}_k^{p-d_{\widetilde{\alpha}}-d_h, r+1}, \mathcal{Z}_\beta),$$
$$\mathcal{Z}_\beta = \{i_0, i_1, \ldots, i_{z_\beta-1}\}.$$

- $\Gamma_2$:

$$\{\phi_{2,j}\}_{j=0}^{n_2-1}: \ (0,0,\widetilde{g}_2) = \left(0,0,N_j^{p-2d_{\widetilde{\alpha}}, r}\right), \quad n_2 = \dim(\mathcal{S}(\mathcal{T}_k^{p-2d_{\widetilde{\alpha}}, r}, [0,1])).$$

**Remark 5.** In the basis construction above we used the fact that suitable bases of $\Gamma_0$ and $\Gamma_1$ are given by the set of linearly independent functions

$$\{N_j^{p,r+2}, \ j = 0,1,\ldots,n_0-1\} \cup \mathcal{N}(\mathcal{T}_k^{p,r+2}, \{1,2,\ldots,k\}) \cup \widetilde{\mathcal{N}}(\mathcal{T}_k^{p,r+2}, \mathcal{Z}_\beta),$$

and

$$\{h\, N_j^{p-d_{\widetilde{\alpha}}-d_h, r+1}, \ j = 0,1,\ldots,n_1-1\} \cup \{h\, N_j, \ N_j \in \mathcal{N}(\mathcal{T}_k^{p-d_{\widetilde{\alpha}}-d_h, r+1}, \mathcal{Z}_\beta)\},$$

respectively.

The selection of some particular values for $z_\beta$ and for the function $h$ can simplify the representation of the basis functions $\phi_{i,j} \in \mathcal{V}_2^2$. This is demonstrated on the basis of the following example.

**Example 2.** Let $q = 1$, which implies that $\widetilde{\alpha}^{(S)} = \alpha^{(S)}$ and $h = 1$. Furthermore, let $z_\beta = 0$, and recall (33) and (34). All except $k$ basis functions have the simple form

$$\{\phi_{0,j}\}_{j=0}^{n_0-1} : g_{0,j}^{(S)}(u,v) = N_j^{p,r+2}(v)M_0(u) + \beta^{(S)}(N_j^{p,r+2})'(v)M_1(u) + \left(\beta^{(S)}(v)\right)^2 (N_j^{p,r+2})''(v)M_2(u),$$
$$\{\phi_{1,j}\}_{j=0}^{n_1-1} : g_{1,j}^{(S)}(u,v) = \widetilde{\alpha}^{(S)}(v)N_j^{p-d_{\widetilde{\alpha}},r+1}(v)M_1(u) + 2\widetilde{\alpha}^{(S)}(v)\beta^{(S)}(v)(N_j^{p-d_{\widetilde{\alpha}},r+1})'(v)M_2(u),$$
$$\{\phi_{2,j}\}_{j=0}^{n_2-1} : g_{2,j}^{(S)}(u,v) = \left(\widetilde{\alpha}^{(S)}(v)\right)^2 N_j^{p-2d_{\widetilde{\alpha}},r}(v)M_2(u),$$

and the remaining $k$ basis functions are given by

$$\begin{aligned}\{\phi_{0,n_0+j}\}_{j=0}^{k-1} : \quad g_{0,j}^{(S)}(u,v) &= N_j(v)M_0(u) + \left(\zeta_{1,j}\widetilde{\alpha}^{(S)}(v) + \beta^{(S)}(v)\right) N_j'(v)M_1(u) \\ &+ \left(\left(\beta^{(S)}(v)\right)^2 + 2\zeta_{1,j}\widetilde{\alpha}^{(S)}(v)\beta^{(S)}(v) + \zeta_{2,j}\left(\widetilde{\alpha}^{(S)}(v)\right)^2\right) N_j''(v)M_2(u),\end{aligned}$$

where $N_j \in \mathcal{N}(\mathcal{T}_k^{p,r+2}, \{1,2,\ldots,k\})$ and $\zeta_{1,j} = -\frac{\beta^{(L)}(\tau_{j+1})}{2\widetilde{\alpha}^{(L)}(\tau_{j+1})} - \frac{\beta^{(R)}(\tau_{j+1})}{2\widetilde{\alpha}^{(R)}(\tau_{j+1})}$, $\zeta_{2,j} = \frac{\beta^{(L)}(\tau_{j+1})\beta^{(R)}(\tau_{j+1})}{\widetilde{\alpha}^{(L)}(\tau_{j+1})\widetilde{\alpha}^{(R)}(\tau_{j+1})}$.

In case of $\boldsymbol{F}^{(L)}, \boldsymbol{F}^{(R)} \in \mathcal{S}(\mathcal{T}_3^{5,2}, [0,1]^2) \times \mathcal{S}(\mathcal{T}_3^{5,2}, [0,1]^2)$ with the uniform knot vector

$$\mathcal{T}_3^{5,2} = (0,0,0,0,0,0,\frac{1}{4},\frac{1}{4},\frac{1}{4},\frac{1}{2},\frac{1}{2},\frac{1}{2},\frac{3}{4},\frac{3}{4},\frac{3}{4},1,1,1,1,1,1),$$

the graphs of the basis functions $\phi_{i,j}$ of the space $\mathcal{V}_2^2$ are shown in Fig. 2.

**Remark 6.** While the basis functions of $\mathcal{V}_1^2$ are independent of the geometry mappings $\boldsymbol{F}^{(L)}$ and $\boldsymbol{F}^{(R)}$, the basis functions $\phi_{i,j}$ of $\mathcal{V}_2^2$ depend on the initial two-patch geometry $\boldsymbol{F}$ via the functions $\alpha^{(L)}$, $\alpha^{(R)}$, $\beta^{(L)}$ and $\beta^{(R)}$. In addition, the values of the spline functions $g_{i,j}^{(S)}$ of the isogeometric functions $\phi_{i,j}$, $i = 1,2$, strongly depend on the degree $p$ and the value of the knot $\tau_1$ via the functions $M_1$ or $M_2$ (compare Fig. 2). More precisely, for a fixed geometry $\boldsymbol{F}$ the absolute values of the spline functions $g_{i,j}^{(S)}$, $i = 1,2$, decrease in case that $p$ increases or the value of $\tau_1$ decreases. Clearly, if desired, an appropriate scaling of these spline functions could adjust this influence.

*Spline coefficients as matrix entries.* Let $\boldsymbol{B}^{(S)} = (\boldsymbol{B}_0^{(S)}, \boldsymbol{B}_1^{(S)}, \boldsymbol{B}_2^{(S)})^T$ be the vector of functions composed of

$$\boldsymbol{B}_0^{(S)}(u,v) = (g_{0,0}^{(S)}(u,v), g_{0,1}^{(S)}(u,v), \ldots, g_{0,n_0+k+z_\beta-1}^{(S)}(u,v))^T,$$

$$\boldsymbol{B}_1^{(S)}(u,v) = (g_{1,0}^{(S)}(u,v), g_{1,1}^{(S)}(u,v), \ldots, g_{1,n_1+z_\beta-1}^{(S)}(u,v))^T$$

and

$$\boldsymbol{B}_2^{(S)}(u,v) = (g_{2,0}^{(S)}(u,v), g_{2,1}^{(S)}(u,v), \ldots, g_{2,n_2-1}^{(S)}(u,v))^T,$$

and let $\boldsymbol{B}^* = (\boldsymbol{B}_0^*, \boldsymbol{B}_1^*, \boldsymbol{B}_2^*)^T$ be the vector of functions composed of

$$\boldsymbol{B}_0^*(u,v) = (N_{0,0}^{p,r}(u,v), N_{0,1}^{p,r}(u,v), \ldots, N_{0,n-1}^{p,r}(u,v))^T,$$

$$\boldsymbol{B}_1^*(u,v) = (N_{1,0}^{p,r}(u,v), N_{1,1}^{p,r}(u,v), \ldots, N_{1,n-1}^{p,r}(u,v))^T$$

and

$$\boldsymbol{B}_2^*(u,v) = (N_{2,0}^{p,r}(u,v), N_{2,1}^{p,r}(u,v), , \ldots, N_{2,n-1}^{p,r}(u,v))^T,$$

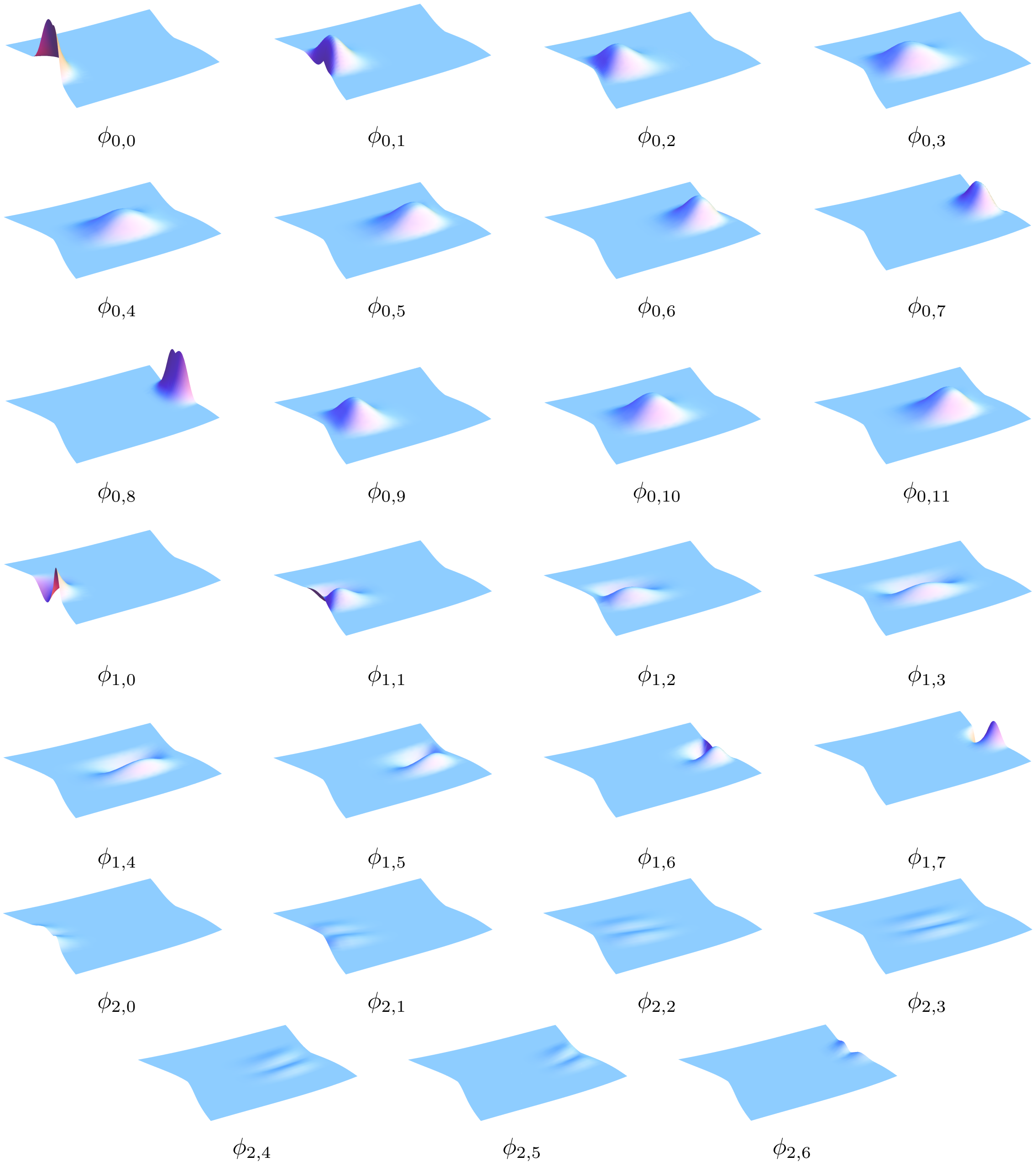

Figure 2: The graphs of the basis functions $\phi_{i,j}$ of the resulting space $\mathcal{V}_2^2$, when both patches $\boldsymbol{F}^{(L)}$ and $\boldsymbol{F}^{(R)}$ are represented in the space $\mathcal{S}(\mathcal{T}_3^{5,2}, [0,1]^2) \times \mathcal{S}(\mathcal{T}_3^{5,2}, [0,1]^2)$ for a uniform knot vector $\mathcal{T}_3^{5,2}$. The graphs are plotted in the parameter range $[0, \frac{1}{2}] \times [0, 1]$ for both patches $\boldsymbol{F}^{(L)}$ and $\boldsymbol{F}^{(R)}$.

where $n = p+1+k(p-r)$. There exists a matrix $A^{(S)} \in \mathbb{R}^{(n_0+n_1+n_2+k+2z_\beta)\times 3n}$ such that

$$\boldsymbol{B}^{(S)}(u,v) = A^{(S)}\boldsymbol{B}^*(u,v). \tag{39}$$

Thereby, the matrix $A^{(S)}$ for $S \in \{L,R\}$ has a block structure, since equation (39) is equivalent to

$$\begin{pmatrix} \boldsymbol{B}_0^{(S)} \\ \boldsymbol{B}_1^{(S)} \\ \boldsymbol{B}_2^{(S)} \end{pmatrix} = \begin{pmatrix} \bar{A}_0 & \widehat{A}_0^{(S)} & \widetilde{A}_0^{(S)} \\ 0 & \widehat{A}_1^{(S)} & \widetilde{A}_1^{(S)} \\ 0 & 0 & \widetilde{A}_2^{(S)} \end{pmatrix} \begin{pmatrix} \boldsymbol{B}_0^* \\ \boldsymbol{B}_1^* \\ \boldsymbol{B}_2^* \end{pmatrix}, \tag{40}$$

where $\bar{A}_0$, $\widehat{A}_0^{(S)}$, $\widetilde{A}_0^{(S)} \in \mathbb{R}^{(n_0+k+z_\beta)\times n}$, $\widehat{A}_1^{(S)}$, $\widetilde{A}_1^{(S)} \in \mathbb{R}^{(n_1+z_\beta)\times n}$ and $\widetilde{A}_2^{(S)} \in \mathbb{R}^{n_2\times n}$. For each row the single entries of the matrix $A^{(S)}$, $S \in \{L,R\}$, are the B-spline coefficients of the corresponding spline function $g_{0,i}^{(S)}$, $g_{1,j}^{(S)}$ or $g_{2,i}^{(S)}$ with respect to the spline space $\mathcal{S}(\mathcal{T}_k^{p,r},[0,1]^2)$.

Similar to [17] for the case of $C^1$-smooth isogeometric basis functions for analysis-suitable $G^1$ two-patch geometries, the single entries of the matrices can be computed by means of the concept of blossoming (cf. [7, 10, 24, 32, 33]) or interpolation. For the sake of simplicity and brevity, we only present the second possibility. The $m$-th rows of the matrices $\bar{A}_0$, $\widehat{A}_{\hat{\imath}}^{(S)}$, $\hat{\imath} = 0,1$, and $\widetilde{A}_{\tilde{\imath}}^{(S)}$, $\tilde{\imath} = 0,1,2$, denoted by $(\bar{a}_{0,m,0},\ldots,\bar{a}_{0,m,n-1})$, $(\widehat{a}_{\hat{\imath},m,0}^{(S)},\ldots,\widehat{a}_{\hat{\imath},m,n-1}^{(S)})$ and $(\widetilde{a}_{\tilde{\imath},m,0}^{(S)},\ldots,\widetilde{a}_{\tilde{\imath},m,n-1}^{(S)})$, respectively, are computed by solving the following systems of linear equations

$$g_{0,m}^{(L)}(0,\xi_i) = \sum_{j=0}^{n-1} \bar{a}_{0,m,j} N_j^{p,r}(\xi_i), \quad i = 0,1,\ldots,n-1,$$

$$\frac{\tau_1 D_u g_{\hat{\imath},m}^{(S)}(0,\xi_i)}{p} + g_{\hat{\imath},m}^{(S)}(0,\xi_i) = \sum_{j=0}^{n-1} \hat{a}_{\hat{\imath},m,j}^{(S)} N_j^{p,r}(\xi_i), \quad i = 0,1,\ldots,n-1,$$

and

$$\frac{\tau_1^2 D_{uu} g_{\tilde{\imath},m}^{(S)}(0,\xi_i)}{p(p-1)} + \frac{2\tau_1 D_u g_{\tilde{\imath},m}^{(S)}(0,\xi_i)}{p} + g_{\tilde{\imath},m}^{(S)}(0,\xi_i) = \sum_{j=0}^{n-1} \widetilde{a}_{\tilde{\imath},m,j}^{(S)} N_j^{p,r}(\xi_i), \quad i = 0,1,\ldots,n-1,$$

respectively, where $(\xi_i)_{i=0,1,\ldots,n-1}$ are the Greville abscissae of the B-splines $N_j^{p,r}$, $j = 0,1,\ldots,n-1$, of the spline space $\mathcal{S}(\mathcal{T}_k^{p,r},[0,1])$. Thereby, these systems of linear equations involve the evaluation of equation (40) at $u=0$, the evaluation of the first derivative of (40) and the evaluation of the second derivative of (40) both with respect to u at $u=0$.

## 7. Approximation properties of the space $\mathcal{V}^2$

We will introduce the subspace $\mathcal{W}^2$ of $\mathcal{V}^2$ for which an even simpler and more uniform basis construction will be presented. Theoretical investigations and numerical experiments will indicate that both spaces $\mathcal{V}^2$ and $\mathcal{W}^2$ enjoy optimal approximation properties.

*The subspace $\mathcal{W}^2$ of $\mathcal{V}^2$.* Recall Corollary 1 and Remark 3. Clearly, any choice of the function triplets

$$(g_0, g_1, g_2) \in \mathcal{S}(\mathcal{T}_k^{p,r+2}, [0,1]) \times \mathcal{S}(\mathcal{T}_k^{p-d_\alpha,r+1}, [0,1]) \times \mathcal{S}(\mathcal{T}_k^{p-2d_\alpha,r}, [0,1]) \qquad (41)$$

in (35) leads to unique functions $g^{(S)} \in \mathcal{S}(\mathcal{T}_k^{p,r}, [0,1]^2)$, $S \in \{L, R\}$, and hence to a unique isogeometric function $\phi \in \mathcal{V}_2^2$. Let us introduce the subspace of $\mathcal{V}_2^2$, which consists of these $C^2$-smooth isogeometric functions $\phi$. We denote by $\mathcal{W}_2^2$ the space

$$\mathcal{W}_2^2 = \{\phi \in \mathcal{V}_2^2 \mid \phi \text{ possesses representations (35) for } g^{(S)} \text{ with function triplets (41)}\},$$

and by $\mathcal{W}^2$ the space

$$\mathcal{W}^2 = \mathcal{V}_1^2 \oplus \mathcal{W}_2^2.$$

Relation $\mathcal{W}_2^2 \subseteq \mathcal{V}_2^2$ implies $\mathcal{W}^2 \subseteq \mathcal{V}^2$. Since any choice (41) leads to a unique isogeometric function $\phi \in \mathcal{W}_2^2$, the dimensions of $\mathcal{W}_2^2$ and $\mathcal{W}^2$ are given as follows:

**Corollary 4.** *The dimensions of $\mathcal{W}_2^2$ and $\mathcal{W}^2$ are equal to*

$$\dim \mathcal{W}_2^2 = (k+1)(3p - 3d_\alpha) + 3(1-k-kr)$$

*and*

$$\dim \mathcal{W}^2 = \underbrace{2(p-2+k(p-r))(p+1+k(p-r))}_{\dim \mathcal{V}_1^2} + \underbrace{(k+1)(3p-3d_\alpha) + 3(1-k-kr)}_{\dim \mathcal{W}_2^2},$$

*respectively.*

Recall that $n_0 = \dim(\mathcal{S}(\mathcal{T}_k^{p,r+2}, [0,1]))$. In addition, let $\bar{n}_1 = \dim(\mathcal{S}(\mathcal{T}_k^{p-d_\alpha,r+1}, [0,1]))$ and $\bar{n}_2 = \dim(\mathcal{S}(\mathcal{T}_k^{p-2d_\alpha,r}, [0,1]))$. Basis functions of the space $\mathcal{W}_2^2$ can be constructed similarly to the case of $\mathcal{V}_2^2$. We choose the function triplets $(g_0, g_1, g_2) = (N_j^{p,r+2}, 0, 0)$, $j = 0, 1, \ldots, n_0 - 1$, $(g_0, g_1, g_2) = (0, N_j^{p-d_\alpha,r+1}, 0)$, $j = 0, 1, \ldots, \bar{n}_1 - 1$, and $(g_0, g_1, g_2) = (0, 0, N_j^{p-2d_\alpha,r})$, $j = 0, 1, \ldots, \bar{n}_2 - 1$, which define via relation (35) a basis of $\mathcal{W}_2^2$. More precisely, this basis is given by the collection of isogeometric functions $\bar{\phi}_{i,j}$, which are determined by

$$\{\bar{\phi}_{0,j}\}_{j=0}^{n_0-1} : g_{0,j}^{(S)}(u,v) = N_j^{p,r+2}(v)M_0(u) + \beta^{(S)}(N_j^{p,r+2})'(v)M_1(u) + \left(\beta^{(S)}(v)\right)^2 (N_j^{p,r+2})''(v)M_2(u),$$
$$\{\bar{\phi}_{1,j}\}_{j=0}^{\bar{n}_1-1} : g_{1,j}^{(S)}(u,v) = \alpha^{(S)}(v)N_j^{p-d_\alpha,r+1}(v)M_1(u) + 2\alpha^{(S)}(v)\beta^{(S)}(v)(N_j^{p-d_\alpha,r+1})'(v)M_2(u),$$
$$\{\bar{\phi}_{2,j}\}_{j=0}^{\bar{n}_2-1} : g_{2,j}^{(S)}(u,v) = \left(\alpha^{(S)}(v)\right)^2 N_j^{p-2d_\alpha,r}(v)M_2(u).$$

Clearly, the approximation properties of the spaces $\mathcal{W}^2$ and $\mathcal{V}^2$ are determined by the ones from the spaces $\mathcal{W}_2^2$ and $\mathcal{V}_2^2$, respectively. Analyzing the functions $g_0$, $g_1$ and $g_2$ of a function $\phi \in \mathcal{W}_2^2$, we observe that $g_0$ is the trace of $\phi$ at the interface $\Gamma$, $g_1$ is a function consisting of first derivatives[1] of $\phi$ at $\Gamma$, and $g_2$ is a function consisting of first and second derivatives of $\phi$ at $\Gamma$. Note that functions $g_0$, $g_1$ and $g_2$ can be independently selected as any spline function from the spaces $\mathcal{S}(\mathcal{T}_k^{p,r+2},[0,1])$, $\mathcal{S}(\mathcal{T}_k^{p-d_\alpha,r+1},[0,1])$ and $\mathcal{S}(\mathcal{T}_k^{p-2d_\alpha,r},[0,1])$, respectively. This indicates that the spaces $\mathcal{W}_2^2$ and $\mathcal{V}_2^2$, and therefore the spaces $\mathcal{W}^2$ and $\mathcal{V}^2$, allow optimal convergence under $h$-refinement when $r \leq p-3$. (In case of $r > p-3$, e.g. the space $\mathcal{S}(\mathcal{T}_k^{p,r+2},[0,1])$ will not be refined since $\mathcal{S}(\mathcal{T}_k^{p,r+2},[0,1]) = \mathcal{P}^p([0,1])$.) A theoretical proof of the approximation properties is clearly an interesting question but deserves further investigation and is beyond the scope of the paper. Instead, the optimal approximation order will be shown numerically.

[1]It was shown in [8], that $g_1$ is the transversal derivative on $\Gamma$ with respect to $\boldsymbol{d} = \boldsymbol{d}^{(L)} = \boldsymbol{d}^{(R)}$ with $\boldsymbol{d}^{(S)} \circ \boldsymbol{F}_0(v) = (D_u \boldsymbol{F}^{(S)}(0,v), \boldsymbol{F}_0'(v))(1, -\beta^{(S)}(v))^T \frac{1}{\alpha^{(S)}(v)}$, $S \in \{L, R\}$.

*$L^2$-approximation.* We perform $L^2$-approximation on two different bilinear-like $G^2$ two-patch parameterizations using the spaces $\mathcal{V}^2$ and $\mathcal{W}^2$.

**Example 3.** We consider the two bilinear-like $G^2$ two-patch geometries introduced in Example 1 and Appendix A, and shown in Fig. A.4(third row). Recall that for both parameterizations $\boldsymbol{F}$ the two geometries mappings $\boldsymbol{F}^{(L)}$ and $\boldsymbol{F}^{(R)}$ are biquintic polynomial patches, i.e. $\boldsymbol{F}^{(L)}, \boldsymbol{F}^{(R)} \in \mathcal{P}^5([0,1]^2) \times \mathcal{P}^5([0,1]^2)$. We generate a sequence of nested spaces $\mathcal{V}^2$ and $\mathcal{W}^2$ with respect to the number of different inner knots $k$ of the underlying spline spaces $\mathcal{S}(\mathcal{T}_k^{5,2},[0,1])$. More precisely, we perform dyadic $h$-refinement by choosing $k = 2^L - 1$, where $L$ is the *level of refinement.* The resulting $C^2$-smooth isogeometric spline spaces $\mathcal{V}^2$ and $\mathcal{W}^2$ for the different levels $L$ are denoted by $\mathcal{V}_h^2$ and $\mathcal{W}_h^2$, respectively, where $h = \mathcal{O}(2^{-L})$. In addition, we denote the corresponding subspaces $\mathcal{V}_1^2$, $\mathcal{V}_2^2$ and $\mathcal{W}_2^2$ by $\mathcal{V}_{1,h}^2$, $\mathcal{V}_{2,h}^2$ and $\mathcal{W}_{2,h}^2$, respectively.

We use the same model problem as in [19, Section 6.2] but with a slightly different notation. Let $\mathcal{V}_h$ be the isogeometric space $\mathcal{V}_h^2$ or $\mathcal{W}_h^2$, and let $I_h = \{0,1,\ldots,\dim \mathcal{V}_h - 1\}$. Furthermore, let $\{\phi_i\}_{i\in I_h}$ be the collection of the basis functions of $\mathcal{V}_h$, and let $g_i^{(S)}$, $S \in \{L,R\}$, be the associated spline functions $\phi_i \circ \boldsymbol{F}^{(S)}$. Given the smooth function

$$f : \Omega \to \mathbb{R}, \quad f(\boldsymbol{x}) = f(x_1,x_2) = 2\cos(2x_1)\sin(2x_2), \tag{42}$$

see Fig. 3, the goal is to approximate $f$ by the function

$$u_h(\boldsymbol{x}) = \sum_{i\in I_h} b_i\phi_i(\boldsymbol{x}), \quad b_i \in \mathbb{R},$$

solving the least squares problem

$$\|u_h - f\|_{L^2}^2 = \int_\Omega (u_h(\boldsymbol{x}) - f(\boldsymbol{x}))^2 \mathrm{d}\boldsymbol{x} \to \min_{b_i,\, i\in I_h}, \tag{43}$$

for the unknown coefficients $\{b_i\}_{i\in I_h}$. The minimization problem (43) is equivalent to the system of linear equations

$$M\boldsymbol{b} = \boldsymbol{f}, \quad \boldsymbol{b} = (b_i)_{i\in I_h},$$

where the entries of the mass matrix $M = (m_{i,j})_{i,j\in I_h}$ and of the load vector $\boldsymbol{f} = (f_i)_{i\in I_h}$ are given by

$$m_{i,j} = \int_\Omega \phi_i(\boldsymbol{x})\phi_j(\boldsymbol{x})\mathrm{d}\boldsymbol{x} = \sum_{S\in\{L,R\}} \int_0^1\int_0^1 g_i^{(S)}(u,v)g_j^{(S)}(u,v)|\det(J^{(S)}(u,v))|\mathrm{d}u\mathrm{d}v,$$

and

$$f_i = \int_\Omega f(\boldsymbol{x})\phi_i(\boldsymbol{x})\mathrm{d}\boldsymbol{x} = \sum_{S\in\{L,R\}} \int_0^1 \int_0^1 f(\boldsymbol{F}^{(S)}(u,v))g_i^{(S)}(u,v)|\det(J^{(S)}(u,v))|\mathrm{d}u\mathrm{d}v,$$

respectively, where $J^{(S)}$, $S \in \{L,R\}$, is the Jacobian of $\boldsymbol{F}^{(S)}$.

Table 2 reports the resulting relative $L^2$-errors with the estimated convergence rates by using the spaces $\mathcal{V}_h^2$ and $\mathcal{W}_h^2$ for approximating the function (42) over the two considered bilinear-like $G^2$ two-patch parameterizations. The numerical results indicate that both spaces possess an optimal aproximation order of $\mathcal{O}(h^6)$. In addition, Table 2 presents the obtained condition numbers of the diagonally scaled mass matrices $M$ (cf. [6]) with the estimated growth rates. For both spaces the used $C^2$-smooth isogeometric basis functions are well-conditioned but with significantly smaller condition numbers for the functions of the subspace $\mathcal{W}_h^2$.

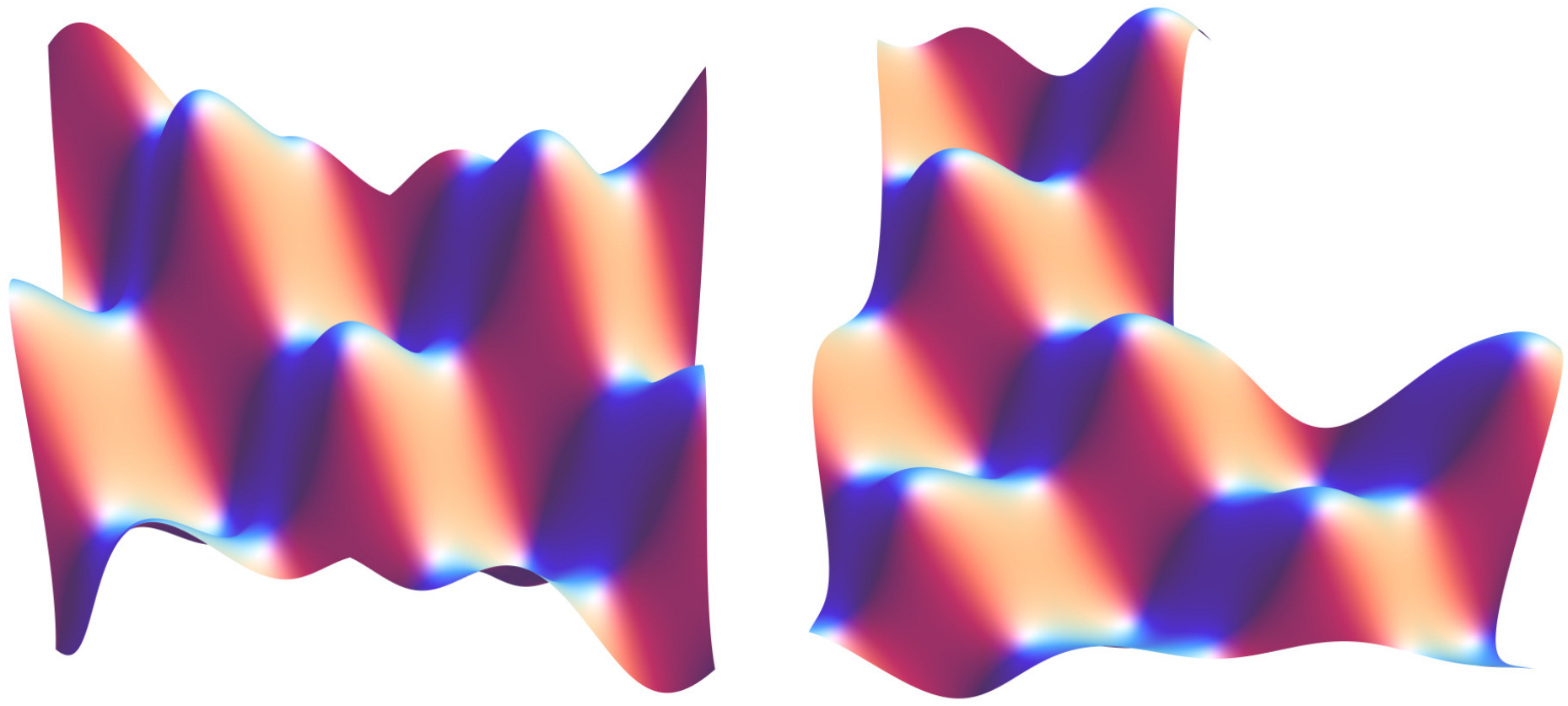

(a) (b)

Figure 3: The function (42) over the two bilinear-like $G^2$ two-patch geometries shown in Fig A.4(third row).

## 8. Conclusion

We studied the space $\mathcal{V}^2$ of $C^2$-smooth isogeometric functions over bilinear-like $G^2$ two-patch parameterizations. A similar approach as presented in [17] for the case of $C^1$-smooth isogeometric functions over the so-called analysis-suitable $G^1$ two-patch parameterizations was generalized here and further developed. A first simple method to generate bilinear-like $G^2$ two-patch parameterizations was presented. Such parameterizations allow curved interfaces and boundaries and examples demonstrate the flexibility of such domains in the sense of the quality of the approximation with respect to given generic two-patch geometries. The dimension of the space $\mathcal{V}^2$ was investigated, and an explicit, locally supported, and well-conditioned basis for the space was presented.

| | | Bilinear-like two-patch geometry (a) | | | | | Bilinear-like two-patch geometry (b) | | | | |
|---|---|---|---|---|---|---|---|---|---|---|---|
| | | Full space $\mathcal{V}^2_{2,h}$ | | | | | | | | | |
| $L$ | $\dim \mathcal{V}^2_{1,h}$ | $\dim \mathcal{V}^2_{2,h}$ | $\frac{\|\|f-u_h\|\|_{L^2}}{\|\|f\|\|_{L^2}}$ | e.c.r. $\|\|\cdot\|\|_{L^2}$ | $\kappa(D^{-\frac{1}{2}}MD^{-\frac{1}{2}})$ | rate $\kappa$ | $\dim \mathcal{V}^2_{2,h}$ | $\frac{\|\|f-u_h\|\|_{L^2}}{\|\|f\|\|_{L^2}}$ | e.c.r. $\|\|\cdot\|\|_{L^2}$ | $\kappa(D^{-\frac{1}{2}}MD^{-\frac{1}{2}})$ | rate $\kappa$ |
| 0 | 36 | 15 | 1.16e-01 | - | 16825.54 | - | 18 | 2.69e-01 | - | 46744.57 | - |
| 1 | 108 | 19 | 7.92e-03 | 3.87 | 32444.61 | -0.95 | 25 | 2.89e-02 | 3.22 | 44746.92 | 0.06 |
| 2 | 360 | 27 | 3.85e-04 | 4.36 | 67575.40 | -1.06 | 39 | 1.47e-03 | 4.30 | 176234.54 | -1.98 |
| 3 | 1296 | 43 | 4.89e-06 | 6.3 | 106706.11 | -0.66 | 67 | 3.59e-05 | 5.35 | 261523.74 | -0.57 |
| 4 | 4896 | 75 | 5.51e-08 | 6.47 | 118077.96 | -0.15 | 123 | 4.68e-07 | 6.26 | 278536.53 | -0.09 |
| 5 | 19008 | 139 | 7.67e-10 | 6.17 | 121572.95 | -0.04 | 235 | 6.25e-09 | 6.26 | 281426.32 | -0.01 |
| | | Subspace $\mathcal{W}^2_{2,h}$ | | | | | | | | | |
| $L$ | $\dim \mathcal{V}^2_{1,h}$ | $\dim \mathcal{W}^2_{2,h}$ | $\frac{\|\|f-u_h\|\|_{L^2}}{\|\|f\|\|_{L^2}}$ | e.c.r. $\|\|\cdot\|\|_{L^2}$ | $\kappa(D^{-\frac{1}{2}}MD^{-\frac{1}{2}})$ | rate $\kappa$ | $\dim \mathcal{W}^2_{2,h}$ | $\frac{\|\|f-u_h\|\|_{L^2}}{\|\|f\|\|_{L^2}}$ | e.c.r. $\|\|\cdot\|\|_{L^2}$ | $\kappa(D^{-\frac{1}{2}}MD^{-\frac{1}{2}})$ | rate $\kappa$ |
| 0 | 36 | 15 | 1.16e-01 | - | 16825.54 | - | 15 | 3.49e-01 | - | 12481.88 | - |
| 1 | 108 | 18 | 8.09e-03 | 3.84 | 32168.00 | -0.93 | 18 | 8.60e-02 | 2.02 | 29913.20 | -1.26 |
| 2 | 360 | 24 | 5.09e-04 | 3.99 | 39914.37 | -0.31 | 24 | 1.78e-02 | 2.28 | 38775.18 | -0.37 |
| 3 | 1296 | 36 | 6.26e-06 | 6.35 | 38809.86 | 0.04 | 36 | 2.14e-04 | 6.37 | 38565.81 | 0.01 |
| 4 | 4896 | 60 | 6.25e-08 | 6.65 | 38083.05 | 0.03 | 60 | 1.18e-06 | 7.49 | 38052.72 | 0.02 |
| 5 | 19008 | 108 | 8.02e-10 | 6.28 | 38006.65 | 0.01 | 108 | 9.83e-09 | 6.92 | 37991.91 | 0.01 |

Table 2: Approximation of the function (42), see Fig. 3, over the two bilinear-like $G^2$ two-patch geometries shown in Fig. A.4(third row) using the spaces $\mathcal{V}^2_h$ and $\mathcal{W}^2_h$. The dimensions of the used spaces $\mathcal{V}^2_{1,h}$, $\mathcal{V}^2_{2,h}$ and $\mathcal{W}^2_{2,h}$, the resulting relative $L^2$-errors with the estimated convergence rates (e.c.r.), and the resulting condition numbers $\kappa$ with the estimated growth rates. Compare Example 3.

Moreover, $L^2$-approximation on different bilinear-like two-patch geometries by using the basis functions of the corresponding spaces $\mathcal{V}^2$ was performed. The obtained numerical results indicate optimal approximation order of the space $\mathcal{V}^2$ and confirm that the generated bases are well-conditioned. Additionally, we introduced a particular subspace $\mathcal{W}^2 \subseteq \mathcal{V}^2$ which already possesses (numerically shown) optimal approximation order, while the basis functions have much simpler representations and smaller condition number.

The paper leaves several open issues which are worth to study. A first interesting topic is the extension of our work to the case of bilinear-like $G^2$ multi-patch parameterizations with extraordinary vertices. Another one is a detailed theoretical investigation of the approximation order of the space of $C^2$-smooth isogeometric functions over bilinear-like two-patch parameterizations. For this purpose, the geometric meaning of the expression $g_2$, similar to the one for $g_0$ and $g_1$ as presented in [8] (also compare Section 7), could be helpful. Some further methods to construct bilinear-like $G^2$ parameterizations (see e.g., [16]) is another interesting issue. Moreover, an extension of the concept of bilinear-like $G^2$ geometries to the perhaps larger set of analysis-suitable $G^2$ geometries should be considered as well. Thereby, the class of analysis-suitable $G^2$ parameterizations should comprises analogous to the $C^1$ case [8] exactly those geometries which allow $C^2$-smooth isogeometric spaces with optimal approximation properties.

Furthermore, an efficient implementation of possible applications of our $C^2$-smooth isogeometric basis functions, such as solving 6-th order partial differential equations by means of isogeometric discretisation is of high importance, too. Possible examples of such partial differential equations are the Phase-field crystal equation [2, 11], the Kirchhoff plate model based on the Mindlin's gradient elasticity theory [29] and the gradient-enhanced continuum damage models [36]. The extension of our approach to volumetric domains and shells where especially locally supported functions are advantageous, should also be considered.

**Acknowledgment.** M. Kapl was partially supported by the European Research Council through the FP7 ERC Consolidator Grant n.616563 HIGEOM, and by the Italian MIUR through the PRIN Metodologie innovative nella modellistica differenziale numerica. V. Vitrih was partially supported by the Slovenian Research Agency (research program P1-0285). These supports are gratefully acknowledged.

## Appendix A. Construction of bilinear-like $G^2$ two-patch parameterizations

We will present the construction of two bilinear-like $G^2$ two-patch parameterizations from given initial (non-bilinear-like) two-patch parameterizations, see Fig. A.4. On the basis of instance (a), we will first describe the used method to generate the bilinear-like $G^2$ two-patch parameterization, and will use it afterwards analogously for instance (b).

*Bilinear-like $G^2$ two-patch parameterization (a).* Let $\widetilde{\boldsymbol{F}}$ be the initial two-patch parameterization shown in Fig. A.4(a), which consists of two geometry mappings $\widetilde{\boldsymbol{F}}^{(L)}, \widetilde{\boldsymbol{F}}^{(R)} \in \mathcal{P}^3([0,1]^2) \times \mathcal{P}^3([0,1]^2)$ of the form

$$\begin{aligned}\widetilde{\boldsymbol{F}}^{(L)}(u,v) \;=\; & \tfrac{1}{150}(75v(2-v-v^2)+u(-450-234v+9v^2+175v^3)+u^2v(-63+261v-\\ & 148v^2)+u^3v(297-845v+498v^2), 450v-u(75+72v-477v^2+280v^3)+\\ & 4u^2(-75+279v-360v^2+206v^3)+u^3(300-919v+963v^2-544v^3))\end{aligned}$$

and

$$\begin{aligned}\widetilde{\boldsymbol{F}}^{(R)}(u,v) \;=\; & (\tfrac{1}{50}(25v(2-v-v^2)+u(175-90v-21v^2+86v^3)-6u^2v(5-19v+14v^2)+\\ & u^3v(-55+182v-127v^2)), \tfrac{1}{200}(600v+6u(-25+33v+21v^2+21v^3)-\\ & 12u^2v(-27+42v+35v^2)+u^3(100-297v+228v^2+369v^3))).\end{aligned}$$

The goal is to construct a bilinear-like $G^2$ two-patch parameterization $\boldsymbol{F}$ consisting of two geometry mappings $\boldsymbol{F}^{(L)}, \boldsymbol{F}^{(R)} \in \mathcal{P}^5([0,1]^2) \times \mathcal{P}^5([0,1]^2)$, which approximates with good quality the initial two patch parameterization $\widetilde{\boldsymbol{F}}$. For this we peform $L^2$-approximation for each coordinate function of the geometry $\widetilde{\boldsymbol{F}}$ by using the basis functions of a particular biquintic $C^2$-smooth isogeometric space $\mathcal{V}^2$. This space is defined with respect to a bilinear two-patch geometry $\widehat{\boldsymbol{F}}$, see Fig. A.4(a), which is a rough approximation of the initial geometry $\widetilde{\boldsymbol{F}}$. In our case we construct the two bilinear patches $\widehat{\boldsymbol{F}}^{(L)}$ and $\widehat{\boldsymbol{F}}^{(R)}$ of the bilinear two-patch parameterization $\widehat{\boldsymbol{F}}$ by interpolating the vertices of the two initial geometry mappings $\widetilde{\boldsymbol{F}}^{(L)}$ and $\widetilde{\boldsymbol{F}}^{(R)}$. This leads to the bilinear two-patch geometry $\widehat{\boldsymbol{F}}$ visualized in Fig. A.4(a), which possesses the bilinear patches

$$\widehat{\boldsymbol{F}}^{(L)}(u,v) = (-\tfrac{u}{3}(9+v), 3v+\tfrac{u}{6}(-3+5v))$$

and

$$\widehat{\boldsymbol{F}}^{(R)}(u,v) = (\tfrac{7u}{2}(1-v)+3uv, -\tfrac{u}{4}(1-v)+3(1-u)v+\tfrac{7uv}{2}),$$

and the linear functions

$$\alpha^{(L)}(v) = -9 - v,\ \alpha^{(R)}(v) = -\tfrac{3}{2}(-7+v),\ \beta^{(L)}(v) = \tfrac{1}{18}(-3+5v),\ \beta^{(R)}(v) = \tfrac{1}{12}(-1+3v). \tag{A.1}$$

By using the space $\mathcal{V}^2$ to approximate both coordinate functions of the initial two-patch parameterization $\widetilde{\boldsymbol{F}}$, we obtain a two-dimensional surface $\boldsymbol{F}$ belonging to the space $\mathcal{V}^2 \times \mathcal{V}^2$. More precisely, the resulting surface $\boldsymbol{F}$ is a bilinear-like $G^2$ two-patch parameterization with the geometry mappings $\boldsymbol{F}^{(L)}$ and $\boldsymbol{F}^{(R)}$ of the form

$$\begin{aligned}\boldsymbol{F}^{(L)}(u,v) \;=\; & \tfrac{1}{16200}(810v(19-10v-8v^2-2v^3+v^4)+9u(-5469-449v+760v^2-372v^3-761v^4+\\ & 71v^5)+u^2(7308-160224v-29650v^2+171723v^3+31872v^4+2191v^5)+\mathcal{O}(u^3),\\ & 162(-1+300v-10v^2+40v^3-45v^4+17v^5)+9u(-738+1578v-460v^2+1140v^3-\\ & 1641v^4+371v^5)+u^2(-39411-8058v+114472v^2-10452v^3-19209v^4+1018v^5)+\mathcal{O}(u^3))^T\end{aligned}$$

and

$$\begin{aligned}\boldsymbol{F}^{(R)}(u,v) \;=\; & \tfrac{1}{7200}(360v(19-10v-8v^2-2v^3+v^4)+6u(3937+313v-1008v^2-212v^3+233v^4+21v^5)+\\ & u^2(3692-94752v+24918v^2+100483v^3-27528v^4+2007v^5)+\mathcal{O}(u^3),\\ & 72(-1+300v-10v^2+40v^3-45v^4+17v^5)+6u(-426+938v-180v^2+540v^3-\\ & 247v^4+201v^5)+u^2(-23819+6894v+70416v^2-47516v^3+10047v^4-126v^5)+\mathcal{O}(u^3))^T,\end{aligned}$$

and the same linear functions $\alpha^{(L)}$, $\alpha^{(R)}$, $\beta^{(L)}$ and $\beta^{(R)}$, as given in (A.1). As one can observe in Fig A.4(c), the obtained geometry $\boldsymbol{F}$ and the initial geometry $\widetilde{\boldsymbol{F}}$ look indistinquishable, and the resulting discrete relative $L^2$-error

$$\epsilon = \frac{\sum_{S\in\{L,R\}}\sum_{i=0}^{10}\sum_{j=0}^{10}(\widetilde{\boldsymbol{F}}^{(S)}(\frac{i}{10},\frac{j}{10}) - \boldsymbol{F}^{(S)}(\frac{i}{10},\frac{j}{10}))^2}{\sum_{S\in\{L,R\}}\sum_{i=0}^{10}\sum_{j=0}^{10}(\widetilde{\boldsymbol{F}}^{(S)}(\frac{i}{10},\frac{j}{10}))^2}$$

with the value $\epsilon = 4.27\text{e-}05$ is very small.

*Bilinear-like $G^2$ two-patch parameterization (b).* It is analogously constructed as described above. The two bicubic patches $\widetilde{\boldsymbol{F}}^{(L)}$ and $\widetilde{\boldsymbol{F}}^{(R)}$ of the initial two-patch parameterization $\widetilde{\boldsymbol{F}}$ shown in Fig. A.4(b) possess the representation

$$\begin{aligned}\widetilde{\boldsymbol{F}}^{(L)}(u,v) \;=\; & (\tfrac{1}{1050}(50(-21+81v-50v^2+32v^3)+u(-1260+513v+3342v^2-3645v^3)+\\ & 15u^2(252-576v+319v^2+75v^3)+u^3(-2520+7227v-6257v^2+1550v^3)),\\ & \tfrac{1}{350}(-350v(-1-5v+3v^2)+6u(315+266v-1022v^2+566v^3)+u^2(770-\\ & 4158v+8001v^2-4013v^3)+u^3(-560+3262v-6349v^2+3347v^3)))\end{aligned}$$

and

$$\begin{aligned}\widetilde{\boldsymbol{F}}^{(R)}(u,v) \;=\; & (\tfrac{1}{1050}(50(-21+81v-50v^2+32v^3)+u(6300-6480v+8256v^2-4926v^3)+\\ & 3u^2(-350+1887v-3235v^2+1698v^3))+u^3(1050-4491v+7099v^2-\\ & 3658v^3),\tfrac{1}{100}(100v(1+5v-3v^2)+3u(-80+392v-716v^2+379v^3)+\\ & u^2(630-2292v+3264v^2-1552v^3)+u^3(-390+1316v-1556v^2+655v^3))),\end{aligned}$$

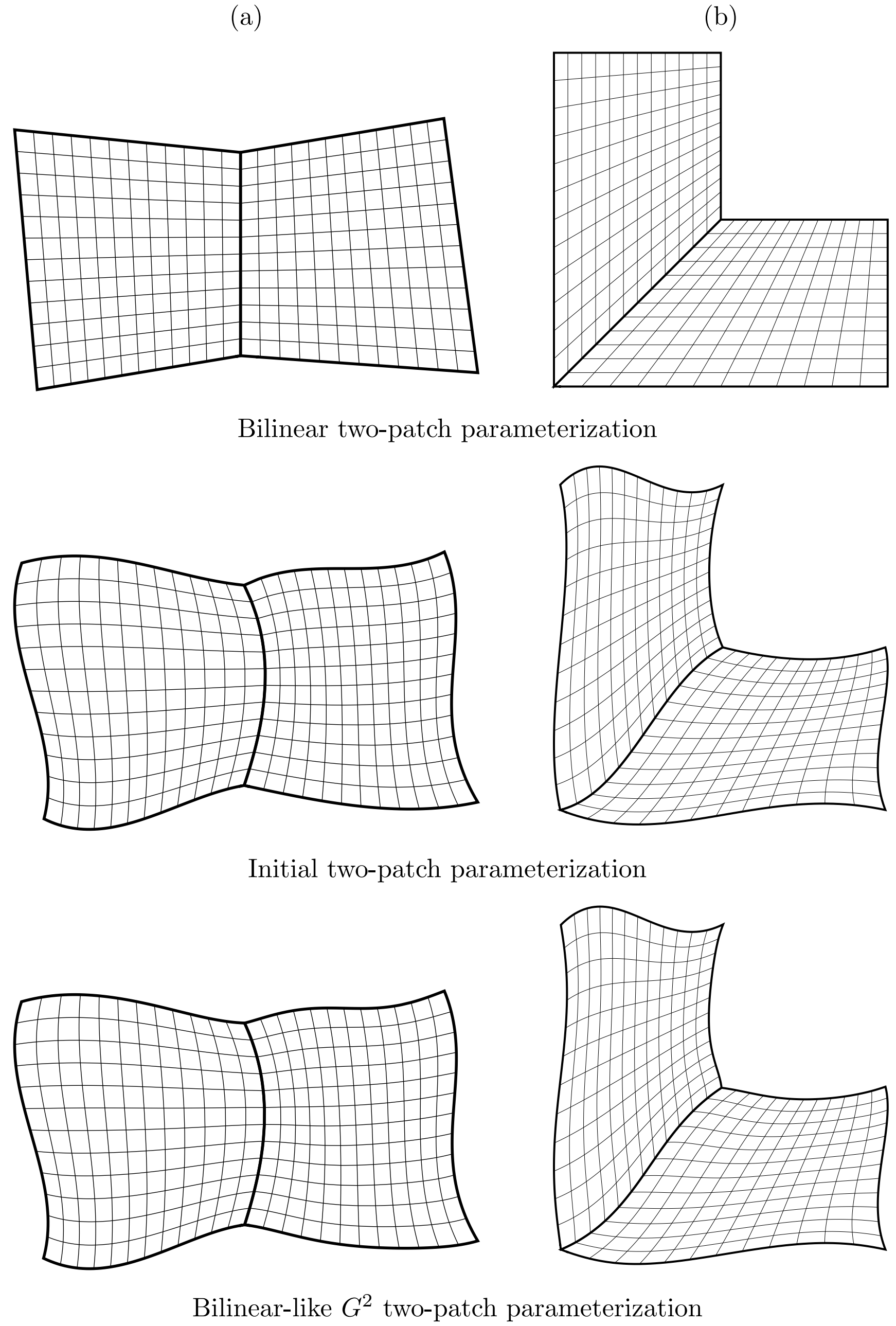

Figure A.4: Construction of bilinear-like $G^2$ two-patch parameterizations $\boldsymbol{F}$ (third row) from initial two-patch parameterizations $\widetilde{\boldsymbol{F}}$ (second row) by means of $L^2$-approximation. The used $C^2$-smooth isogeometric space $\mathcal{V}^2$ in the fitting process is defined on a bilinear reference geometry $\widehat{\boldsymbol{F}}$ (first row).

respectively. The bilinear reference geometry $\widehat{\boldsymbol{F}}$ visualized in Fig. A.4(b) is determined by

$$\widehat{\boldsymbol{F}}^{(L)}(u,v) = (-1+3v, -3u(-2+v)+3v), \quad \widehat{\boldsymbol{F}}^{(R)}(u,v) = (-1-3u(-2+v)+3v, 3v),$$

and is used again to define the biquintic $C^2$-smooth isogeometric space $\mathcal{V}^2$ for the $L^2$-approximation. Again, the resulting bilinear-like $G^2$ two-patch parameterization $\boldsymbol{F}$, see Fig. A.4(b), possesses the same linear functions $\alpha^{(L)}, \alpha^{(R)}$, $\beta^{(L)}$ and $\beta^{(R)}$ as the bilinear reference geometry $\widehat{\boldsymbol{F}}$, given by

$$\alpha^{(L)}(v) = 9(-2+v), \quad \alpha^{(R)}(v) = -9(-2+v), \quad \beta^{(L)}(v) = \beta^{(R)}(v) = 1 - \tfrac{v}{2}.$$

Moreover, the two geometry mappings $\boldsymbol{F}^{(L)}, \boldsymbol{F}^{(R)} \in \mathcal{P}^5([0,1]^2) \times \mathcal{P}^5([0,1]^2)$ are of the form

$$\begin{aligned}\boldsymbol{F}^{(L)} \;=\; & \tfrac{1}{200}(-200+610v+540v^2-2080v^3+2420v^4-896v^5+u(-128+1675v-8220v^2+14780v^3- \\ & 10940v^4+2780v^5)+u^2(-306-6420v+21885v^2-6770v^3-30880v^4+21709v^5)+\mathcal{O}(u^3), \\ & -2+260v+640v^2+140v^3-680v^4+238v^5+u(1034+340v+680v^2-4190v^3+ \\ & 3630v^4-973v^5)+u^2(703-445v-11810v^2+10525v^3+5240v^4-4562v^5)+\mathcal{O}(u^3))^T\end{aligned}$$

and

$$\begin{aligned}\boldsymbol{F}^{(R)} \;=\; & \tfrac{1}{200}(-200+610v+540v^2-2080v^3+2420v^4-896v^5+u(1348-125v-5340v^2+10820v^3- \\ & 7700v^4+1700v^5)+u^2(-1368-660v+8565v^2+10150v^3-41140v^4+23869v^5)+\mathcal{O}(u^3), \\ & -2+260v+640v^2+140v^3-680v^4+238v^5-u(514-1960v+1120v^2+1670v^3- \\ & 1470v^4+217v^5)-u^2(-1549+4045v+3350v^2+635v^3-12260v^4+6074v^5)+\mathcal{O}(u^3))^T,\end{aligned}$$

and the resulting relative $L^2$-error $\epsilon$ is equal to 2.37e-05.

If the approximation quality is not good enough one could either increase the degree of the approximant or use instead of polynomials $C^2$-smooth splines with inner knots. Instead of $L^2$ approximation one could also try to interpolate some boundary data of the initial domain parameterization in order to define bilinear-like $G^2$ parameterization (see e.g. [16]).

## References

[1] F. Auricchio, L. Beirão da Veiga, A. Buffa, C. Lovadina, A. Reali, and G. Sangalli. A fully "locking-free" isogeometric approach for plane linear elasticity problems: a stream function formulation. *Computer Methods in Applied Mechanics and Engineering*, 197(1):160–172, 2007.

[2] A. Bartezzaghi, L. Dedè, and A. Quarteroni. Isogeometric analysis of high order partial differential equations on surfaces. *Computer Methods in Applied Mechanics and Engineering*, 295:446 – 469, 2015.

[3] L. Beirão da Veiga, A. Buffa, G. Sangalli, and R. Vázquez. Mathematical analysis of variational isogeometric methods. *Acta Numerica*, 23:157–287, 5 2014.

[4] D. J. Benson, Y. Bazilevs, M.-C. Hsu, and T. J.R. Hughes. A large deformation, rotation-free, isogeometric shell. *Computer Methods in Applied Mechanics and Engineering*, 200(13):1367–1378, 2011.

[5] M. Bercovier and T. Matskewich. Smooth Bézier surfaces over arbitrary quadrilateral meshes. Technical Report 1412.1125, arXiv.org, 2014.

[6] A. M. Bruaset. *A survey of preconditioned iterative methods*, volume 328 of *Pitman Research Notes in Mathematics Series*. Longman Scientific & Technical, Harlow, 1995.

[7] X. Chen, R. F. Riesenfeld, and E. Cohen. An algorithm for direct B-spline multiplication. *IEEE Transactions on Automation Science and Engineering*, 6(3):433 – 442, 2009.

[8] A. Collin, G. Sangalli, and T. Takacs. Analysis-suitable $G^1$ multi-patch parametrizations for $C^1$ isogeometric spaces. *Computer Aided Geometric Design*, 47:93 – 113, 2016.

[9] J. A. Cottrell, T.J.R. Hughes, and Y. Bazilevs. *Isogeometric Analysis: Toward Integration of CAD and FEA*. John Wiley & Sons, Chichester, England, 2009.

[10] R. Goldman. *Pyramid algorithms : a dynamic programming approach to curves and surfaces for geometric modeling*. Morgan Kaufmann, San Francisco (Calif.), 2003.

[11] H. Gomez and X. Nogueira. An unconditionally energy-stable method for the phase field crystal equation. *Computer Methods in Applied Mechanics and Engineering*, 249 – 252:52 – 61, 2012.

[12] D. Groisser and J. Peters. Matched $G^k$-constructions always yield $C^k$-continuous isogeometric elements. *Computer Aided Geometric Design*, 34:67 – 72, 2015.

[13] J. Hoschek and D. Lasser. *Fundamentals of computer aided geometric design*. A K Peters Ltd., Wellesley, MA, 1993.

[14] T. J. R. Hughes, J. A. Cottrell, and Y. Bazilevs. Isogeometric analysis: CAD, finite elements, NURBS, exact geometry and mesh refinement. *Computer Methods in Applied Mechanics and Engineering*, 194(39-41):4135–4195, 2005.

[15] M. Kapl, F. Buchegger, B. Bercovier, and B. Jüttler. Isogeometric analysis with geometrically continuous functions on planar multi-patch geometries. *Computer Methods in Applied Mechanics and Engineering*, 316:209 – 234, 2017.

[16] M. Kapl, G. Sangalli, and T. Takacs. Construction of analysis-suitable $G^1$ planar multi-patch parameterizations. Technical Report 1706.03264, arXiv.org, 2017.

[17] M. Kapl, G. Sangalli, and T. Takacs. Dimension and basis construction for analysis-suitable $G^1$ two-patch parameterizations. *Computer Aided Geometric Design*, 52–53:75–89, 2017.

[18] M. Kapl and V. Vitrih. Space of $C^2$-smooth geometrically continuous isogeometric functions on planar multi-patch geometries: Dimension and numerical experiments. *Comput. Math. Appl.*, 73(10):2319–2338, 2017.

[19] M. Kapl and V. Vitrih. Space of $C^2$-smooth geometrically continuous isogeometric functions on two-patch geometries. *Comput. Math. Appl.*, 73(1):37 – 59, 2017.

[20] M. Kapl, V. Vitrih, B. Jüttler, and K. Birner. Isogeometric analysis with geometrically continuous functions on two-patch geometries. *Comput. Math. Appl.*, 70(7):1518 – 1538, 2015.

[21] J. Kiendl, Y. Bazilevs, M.-C. Hsu, R. Wüchner, and K.-U. Bletzinger. The bending strip method for isogeometric analysis of Kirchhoff-Love shell structures comprised of multiple patches. *Computer Methods in Applied Mechanics and Engineering*, 199(35):2403–2416, 2010.

[22] J. Kiendl, K.-U. Bletzinger, J. Linhard, and R. Wüchner. Isogeometric shell analysis with Kirchhoff-Love elements. *Computer Methods in Applied Mechanics and Engineering*, 198(49):3902–3914, 2009.

[23] M.-J. Lai and L. L. Schumaker. *Spline functions on triangulations*, volume 110 of *Encyclopedia of Mathematics and its Applications.* Cambridge University Press, Cambridge, 2007.

[24] W. Liu. A simple, efficient degree raising algorithm for B-spline curves. *Computer Aided Geometric Design*, 14(7):693 – 698, 1997.

[25] B. Mourrain, R. Vidunas, and N. Villamizar. Dimension and bases for geometrically continuous splines on surfaces of arbitrary topology. *Computer Aided Geometric Design*, 45:108 – 133, 2016.

[26] T. Nguyen, K. Karčiauskas, and J. Peters. A comparative study of several classical, discrete differential and isogeometric methods for solving poisson's equation on the disk. *Axioms*, 3(2):280–299, 2014.

[27] T. Nguyen, K. Karčiauskas, and J. Peters. $C^1$ finite elements on non-tensor-product 2d and 3d manifolds. *Applied Mathematics and Computation*, 272:148 – 158, 2016.

[28] T. Nguyen and J. Peters. Refinable $C^1$ spline elements for irregular quad layout. *Computer Aided Geometric Design*, 43:123 – 130, 2016.

[29] J. Niiranen, J. Kiendl, A. H. Niemi, and A. Reali. Isogeometric analysis for sixth-order boundary value problems of gradient-elastic Kirchhoff plates. *Computer Methods in Applied Mechanics and Engineering*, 316:328–348, 2017.

[30] J. Peters. Smooth mesh interpolation with cubic patches. *Computer-Aided Design*, 22(2):109 – 120, 1990.

[31] J. Peters. Geometric continuity. In *Handbook of computer aided geometric design*, pages 193–227. North-Holland, Amsterdam, 2002.

[32] L. Ramshaw. Blossoms are polar forms. *Computer Aided Geometric Design*, 6(4):323–358, 1989.

[33] H.-P. Seidel. An introduction to polar forms. *IEEE Computer Graphics and Applications*, 13(1):38–46, 1993.

[34] A. Tagliabue, L. Dedè, and A. Quarteroni. Isogeometric analysis and error estimates for high order partial differential equations in fluid dynamics. *Computers & Fluids*, 102:277 – 303, 2014.

[35] D. Toshniwal, H. Speleers, R. R. Hiemstra, and T. J. R. Hughes. Multi-degree smooth polar splines: A framework for geometric modeling and isogeometric analysis. *Computer Methods in Applied Mechanics and Engineering*, 316:1005–1061, 2017.

[36] C. V. Verhoosel, M. A. Scott, T. J. R. Hughes, and R. de Borst. An isogeometric analysis approach to gradient damage models. *Internat. J. Numer. Methods Engrg.*, 86(1):115–134, 2011.